\documentclass[11pt]{amsart}

\renewcommand{\bar}{\overline}

\newcommand{\pa}{\partial}
\renewcommand{\phi}{\varphi}
\newcommand{\mi}{{\mathcal M}}

\newcounter{hours}\newcounter{minutes}

\newcommand{\ka}{K\"ahler }

\newcommand{\C}{{\mathbb C}}
\newcommand{\Z}{{\mathbb Z}}
\newcommand{\R}{{\mathbb R}}

\newcommand{\frk}[1]{{\mathfrak{#1}}}

\newcommand{\bb}{{\frac{\sqrt{-1}}{2\pi}}}

\usepackage{amsmath,amsthm,amscd} 
\usepackage{calc}        

\title 
[]{Weil-Petersson geometry on moduli space of polarized Calabi-Yau
manifolds}
\author{Zhiqin Lu and Xiaofeng Sun}
\date{April 25, 2002}

\address[Zhiqin Lu and Xiaofeng Sun] 
{Department of Mathematics\\
University of California, Irvine\\ 
Irvine, CA 92697}
\email[Zhiqin Lu]{zlu@math.uci.edu}
\email[Xiaofeng Sun]{xsun@math.uci.edu}
\thanks{The first author is supported by NSF grant
DMS 0204667 and an Alfred P. Sloan fellowship. 
The second author
is supported by NSF grant DMS 0202508.}

\newtheorem{theorem}{Theorem}[section] 

\newtheorem{lemma}{Lemma}[section]
\newtheorem{cor}{Corollary}[section]
\newtheorem{prop}{Proposition}[section]

\newtheorem{definition}{Definition}[section]

\theoremstyle{remark}
\newtheorem{rem}{Remark}[section]     

\begin{document}   
\maketitle


\numberwithin{equation}{section}         

\tableofcontents

\section{Introduction}
Moduli spaces of general polarized algebraic varieties are
studied extensively by algebraic geometers. 
However, there are two classes of moduli spaces where the methods
of differential geometry are equally powerful. These are the
moduli spaces of curves and the moduli spaces of polarized
Calabi-Yau manifolds. Both spaces are complex orbifolds. The
Weil-Petersson metric is the main tool for investigating the
geometry of such moduli spaces.  Under the Weil-Petersson
metrics, these moduli spaces are  \ka orbifolds.

The GIT construction of the (coarse) moduli space (see ~\cite{V2})
of Mumford is  as follows:
let $X$ be a Calabi-Yau manifold and let $L$ be an
ample line bundle over $X$. The pair $(X,L)$ is called
a polarized Calabi-Yau manifold.  Choose a large $m$
such that $L^m$ is very ample. In this way $X$ is
embedded into a complex projecive space $\mathbb C
P^N$. Let
 $\frak{Hilb}(X)$  be the Hilbert scheme of $X$. It is
a compact complex variety. The group $G=PSL(N+1,C)$
acts on  $\frak{Hilb}(X)$ and the moduli space
$\mathcal M$ is the quotient of the stable points of 
$\frak{Hilb}(X)$
by the group $G$. For the purpose of this paper, we
assume that $\mathcal M$ is connected.

The curvature of these moduli spaces with respect to the
Weil-Petersson metric has been studied by many people. For the moduli
space of curves, Wolpert~\cite{Wolpert} gave an explicit formula
for the curvature and proved that the (Riemannian) sectional curvature
 of the Weil-Petersson is negative. Siu~\cite{Siu2} generalized the
result to  the moduli spaces of K\"ahler-Einstein manifolds with
$c_1<0$. Schumacher~\cite{Sc1}, using Siu's
methods, computed the curvature tensor of the moduli spaces of
K\"ahler-Einstein manifolds in the case of $c_1>0$ and $c_1< 0$
respectively.\footnote{Schumacher's method also yields
the curvature formula in the case of $c_1=0$.} Furthermore, 
Strominger~\cite{S} gave the curvature formula for the moduli space
of Calabi-Yau threefolds using the Yukawa couplings. 
Generalizing the formula, C.
Wang~\cite{Wang1} proved the curvarure formula for Calabi-Yau
n-folds where there are no Yukawa couplings. His proof is
purely Hodge theoretic  and is also true on Weil-Petersson
varieties.

It is
important and interesting to know the geometry of  moduli
space at infinity. In~\cite{JY1}, Jost and Yau were able to
understand  the moduli spaces of curves at infinity using the
Schwarz-Yau lemma~\cite{Y3}.   For  moduli space of
polarized Calabi-Yau manifolds, similar results could  be found
in~\cite{LTYZ}. In order to  make use of the Schwarz-Yau lemma, we
need  some natural metric  on the moduli spaces whose holomorphic
sectional  curvature is negative away from zero.

Unlike the 
case of moduli space of curves, the sectional curvature of the
Weil-Petersson metric on moduli space of polarized Calabi-Yau
manifolds is not negative, even in the case when the moduli space is
one dimensional. The curvature of the Weil-Petersson metric can 
either be positive or negative (cf. ~\cite[page 65]{COGP}) on the
moduli space of Calabi-Yau threefolds which are mirror manifolds of
the quintic hypersurfaces in
$CP^4$. This fact prevents us from using the Schwarz-Yau lemma
directly.
 
In~\cite{Lu3}, the first author introduced the Hodge
metric
on the moduli space of polarized Calabi-Yau manifolds.
The Hodge metric is a \ka metric on the moduli space. Its
holomorphic bisectional curvature is nonpositive and both of its
Ricci and holomorphic sectional  curvature are negative away from
zero.  The Hodge metric on moduli space of
Calabi-Yau manifolds is the   counterpart of the Weil-Petersson
metric on  Teichm\"uller space. In \S \ref{3half}, 
we took a further step by defining the ``partial Hodge metric''.
We computed the
curvature of the ``partial Hodge metric''. 
The formula is parallel to the curvature formula of
Wolpert~\cite{Wolpert} on Teichm\"uller space. In the case of the
moduli space of Calabi-Yau threefolds and fourfolds, 
we proved that
the ``partial
Hodge metric'' is the same as the Hodge metric, up to a constant.

Perhaps it is useful to make further comments on the motivations of
this paper. 
We go
back to the idea of Griffiths. In~\cite{Gr1},~\cite{Gr2}, Griffiths
defined the period map. It is a holomorphic map from a
moduli space to the ``classifying space'' defined by Griffiths. The
image of the period map is an integral subvariety of the horizontal
distribution by the Griffiths transversality.  The idea of Griffiths
is that by studying the integral submanifold of the horizontal
distribution, one can partially recover the properties of the moduli
space without  having the knowledge of the varieties the moduli
space  parameterized. 

In the case of moduli space of polarized Calabi-Yau manifolds, we
can 
do better. By a theorem of Tian~\cite{T1}, the Weil-Petersson metric
can be defined by the curvature of the first Hodge bundle. This
implies that the Weil-Petersson metric can be defined without the
detailed 
knowledge of  the Calabi-Yau manifolds. The presence of the
Weil-Petersson
metric gives severe  restrictions on integral submanifold of the 
horizontal distribution.

In \S ~\ref{WPG}, we define the Weil-Petersson geometry. This is
defined to be an integral submanifold of the horizontal distribution
with  the Weil-Petersson metric on it. We further 
axiomatize the results of Viehweg~\cite{V2} and Schmid~\cite{Schmid}
in defining the Weil-Petersson geometry. 
Of course, the axioms will give further restrictions of the integral
submanifolds of the horizontal distribution. It has not been
comprehensively studied how these results interact with the geometry
of the integral submanifolds with the Weil-Petersson metrics. 

One of 
the motivation of this paper is to make a firm foundation to  study
these interactions.

Before giving the main results of this paper, we give a short
definition of  Weil-Petersson, Hodge, and partial Hodge metrics.
For detailed definitions, see~\cite{T1, To} for  Weil-Petersson
metric, ~\cite{Lu3, Lu5} for  Hodge metric, and \S 4 for 
 partial  Hodge metric.

All of these three metrics are Hodge theoretic in the sense
that they depend on the variation of the Hodge structures only.
Let $F^n$ be the first Hodge bundle over $\mathcal M$. Then
the Weil-Petersson metric is defined as
\[
\omega_{_{WP}}=c_1(F^n)=-\frac{\sqrt{-1}}{2\pi}\pa\bar\pa\log
Q(\Omega,\bar\Omega),
\]
where $\Omega$ is a local holomorphic nonzero 
section of $F^n$ and $Q$
is the  polariztion (see ~\eqref{polar1}). The Hodge metric is defined
as follows: given the period map
\[
\mathcal M\rightarrow D,
\]
where $D$ is the Griffiths' classifying space. Let
$D=G/V$ for the real semisimple group $G$ defined by the
polarization $Q$. Let $K$ be the connected components of the
maximal compact subgroup of $G$ containing $V$. The space $G/K$ is a 
(Riemannian) symmetric space which carries the unique invariant
metric $ds^2$ (up to a constant). Let $\pi$ be the composition map
$\mathcal M\rightarrow D\rightarrow G/K$. Then the Hodge metric
is defined by $\pi^*(ds^2)$. It is K\"ahlerian.

The partial Hodge metric is defined by
\[
\omega_\mu=\mu\omega_{_{WP}}+{\rm Ric} (\omega_{_{WP}})
\]
for positive number $\mu>m+1$, where $m$ is the dimension of the
moduli space. In the case of Calabi-Yau three and four-folds,
with the suitable choice of $\mu$, partial Hodge metric is 
the Hodge metric. If the dimension of the Calabi-Yau manifolds
is greater than or equal to five, there is no direct link between
the Hodge and the partial Hodge meric.

As the first result of this paper,  
we have the following 
explicit formula for the curvature of the partial Hodge metric:

\begin{theorem}\label{curvmu1}
Let $\mi$ be a moduli space of polarized
Calabi-Yau manifolds. Let the dimension of
$\mi$ be $m$. Let $\omega_{WP}$ be the
\ka form of the Weil-Petersson metric. 
Then the metric $\omega_{\mu}=\mu\omega_{WP}+{\rm Ric}(\omega_{WP})$
is
\ka for $\mu>m+1$ and the curvature tensor of 
$\omega_{\mu}$ is 
\begin{align}\label{gencurv1}
\begin{split}
\widetilde{R}_{i\bar j k\bar l} =& (\mu -m-1)(g_{i\bar j}g_{k\bar l}+
g_{i\bar l}g_{k\bar j})-(\mu -m)F_{i\bar j k\bar l}
+F_{i\bar q \alpha\bar l}F_{p\bar j k\bar\beta}
g^{\alpha\bar\beta}g^{p\bar q}\\
&+F_{\alpha\bar q k\bar l}F_{i\bar j p\bar\beta}
g^{\alpha\bar\beta}g^{p\bar q}
+\frac{(D_{k}D_{\alpha}D_{i}\Omega,
\bar{D_{l}D_{\beta}D_{j}\Omega})}{(\Omega,\bar\Omega)}g^{\alpha\bar\beta}\\
&+\frac{(E_{k\alpha i},
\bar{E_{l\beta j}})}{(\Omega,\bar\Omega)}g^{\alpha\bar\beta}
-h^{s\bar t}\frac{(E_{k\alpha i},\bar{D_{\beta}D_{t}\Omega})}
{(\Omega,\bar\Omega)}\frac{(D_{\gamma}D_{s}\Omega,
\bar{E_{l\tau j}})}{(\Omega,\bar\Omega)}g^{\alpha\bar\beta}
g^{\gamma\bar\tau}.\\
\end{split}
\end{align}
\end{theorem}
(For notations, see \S ~\ref{3half}).

The obvious feature of the above expression is
that
the  high order terms of $R_{i\bar ii\bar i}$ dominates
the high order terms of the rest of the curvature tensor.
Using this, we
can control the Riemannian sectional curvature by the scalar
curvature in the case of Calabi-Yau threefolds (cf. ~\cite{Lu5}) and
of Calabi-Yau fourfolds (Theorem ~\ref{uplow}).

In the case of moduli space of Calabi-Yau fourfolds,  we have the
following result in
\S~\ref{3half}:

\begin{theorem}\label{thm1112}
We use
 the notations as in the above theorem. Let $\mu=m+2$. Then the
bisectional curvature of the \ka metric $\omega_\mu$ is nonpositive.
The Ricci and the holomorphic sectional curvature are all negatively
bounded by the constant $-\frac{1}{m+4}$, where $m$ is the
complex dimension of the moduli space. Furthermore, the partial Hodge
metric is the Hodge metric in the  case of moduli space of Calabi-Yau
fourfolds, up to a constant.
\end{theorem}

\begin{rem}
The Hodge metric was first defined in ~\cite{Lu3}. Using 
Theorem~\ref{hodge20}, one can prove that it is K\"ahler. The fact that 
the holomorphic sectional curvature is negative away from
zero also follows from the classical paper of Griffiths and 
Schmid~\cite{GS}. The nonpositivity of the holomorphic bisectional
curvature is  from ~\cite{Lu3}. The contribution here is that
we find the explicit relation between the Hodge metric and the 
Weil-Petersson metric in the moduli space of Calabi-Yau fourfolds,
and 
 we find out the optimal constant for the upper bound of the 
holomorphic sectional curvature of the Hodge metric.
\end{rem}

We remark that the corresponding result of 
Theorem~\ref{thm1112} in the case of  Calabi-Yau
threefold was proved in~\cite{Lu5}. In the fourfold case, we
don't have the result of Bryant-Griffiths~\cite{BG} about the
integral submanifold of the  horizontal distribution. However, we 
are
still able to prove that in the case of fourfold, the ``partial Hodge
metric'' is the Hodge metric.

Using the above theorems and the Schwarz-Yau lemma, we have the
following global result in \S~\ref{S5}:

\begin{theorem}
Let ${\mathcal M}$ be the moduli space of the polarized Calabi-Yau
manifolds. Then the Hodge volume on any subvarieties of
${\mathcal M}$ is finite. The Riemannian sectional curvature is
$L^1$ bounded with respect to the Weil-Petersson metric on any
subvarieties. In particular, the moduli space of the
polarized Calabi-Yau manifolds has finite Weil-Petersson volume. 
\end{theorem}

\begin{rem}
Since we don't know the boundedness of the curvatures of the Weil-
Petersson or Hodge metric at infinity, it seems to be interesting
to prove that the integral of the curvature is bounded. 
In the one dimensional case, if a complete Riemann surface
has bounded total Gauss curvature, then it is $S^2$
removing finite many points. In high dimensions, we wish to 
find the geometric implications of the fact that the total
curvature is finite.

A more ambitious problem is to prove that the volume
and the integration of the curvatures of the Weil-Petersson 
metric are rational numbers. The same problem on the moduli
space of curves was studied by many people (cf. ~\cite{KMZ},
~\cite{MZograf},~\cite{Sc2}, ~\cite{Zograf1}, ~\cite{Zograf2} 
and ~\cite{Zograf3}).
The difficulty in the case of  moduli space of Calabi-Yau manifolds is
that the compactification is not known to be ``good'' in the sense of
Mumford~\cite[Section 1]{Mumford1}.
 The results of the Weil-Petersson volume on moduli space of
Calabi-Yau manifolds will  be in our next paper~\cite{LS-2}.
\end{rem} 

In the second part of this paper, we study the asymptotic behavior
of the curvature of the Hodge metric at infinity for  moduli
space of dimension one. The problem is related to the
compactification of the moduli space of Calabi-Yau manifolds. By the
theorem of Viehweg~\cite{V2}, the moduli space is a quasi-projective
variety. Other than this result, we don't know much of the asymptotic
behavior  of the moduli space. S. T. Yau suggested that one can
compactify the moduli space by
completing the moduli space using the Weil-Petersson metric first and
then compactifying it. Under his suggestion, we
study the problem. It seems to us that it is
easier to complete the moduli space using
the Hodge metric. After the completion of the
moduli space  using the Hodge metric, one 
would get a metric space which is not worse
than a complex  orbifold. We wish to study  the
curvature of the Hodge metric  near the
infinity of the moduli space in order to study
the Siegel-type theorem~\cite{Mok1} and wish, by using this, we
can give a differential geometric proof of the compactification
theorem of Viehweg. The full results will appear 
at~\cite{LS-2}. In this paper, we have the
following

\begin{theorem}\label{bound-1}
Assume the moduli space $\mathcal M$ of polarized 
Calabi-Yau threefolds is one dimensional. 
If $\Delta^*$ is a holomorphic chart of $\mathcal M$ such 
that $\Delta^*$ 
is complete at $0$ with respect to the Hodge
metric, then the Gauss curvature of the Hodge metric is
bounded.\footnote{The referee pointed out that the result is also true
for partial Hodge metric.} 
\end{theorem}

{\bf Acknowledgment.} We thank Professor P. Li, R.  Schoen and
G. Tian for the encouragement during the preparation of this
paper. The first author especially thanks  Professor Phong for the
interest of the work and the support. Both authors thank Professor
K. Liu for his interest in our work and the discussions.
Finally, we thank the referee for the extremely valuable comments
and criticisms on the manuscript. 
Without his (her) suggestions, the paper won't be in its current
form. 
\section{Preliminaries}\label{2}
Let $X$ be a compact K\"ahler manifold
of dimension $n$. A $C^\infty$ form on $X$ 
decomposes into $(p,q)$-components according to the number of 
$dz's$ and $d\bar{z}'s$. Denoting the $C^\infty$ $n$-forms and
the 
$C^\infty(p,q)$ forms on $X$ by $A^n(X)$ and $A^{p,q}(X)$ 
respectively, we have the following decomposition
\[
A^n(X)=\underset{p+q=n}{\oplus}A^{p,q}(X).
\]
The cohomology group is defined as
\begin{align*}
H^{p,q}(X)=&\{closed \,(p,q) -forms\}/\{exact\, (p,q)-forms\}\\
=&\{\phi\in A^{p,q}(X)|d\phi=0\}/dA^{n-1}(X)\cap A^{p,q}(X).
\end{align*}

The relations between the groups $\{H^{p,q}(X)\}$ and the de Rham
cohomology is the following Hodge decomposition:

\begin{theorem}  (Hodge decomposition theorem) Let $X$ be a 
compact K\"ahler manifold of dimension $n$. Then the n-th 
complex de Rham cohomology group of $X$ can be written as the
direct sum
\begin{equation}\label{hdcomp}
H^n(X,{\Z})\otimes {\C}    
=H^n_{DR}(X,\C)=
\underset{p+q=n}{\oplus}H^{p,q}(X).
\end{equation}
\end{theorem}

A $(1,1)$ form  $\omega$ is called a polarization of $X$ 
if $[\omega]$ is the first Chern class of an ample line bundle over $X$.
The pair 
 $(X,\omega)$ is
called a  polarized algebraic variety.

Using $\omega$, one can define
\[
L: H^k(X,\C)\rightarrow H^{k+2}(X,\C), 
\quad
[\alpha]\mapsto [\alpha\wedge\omega]
\]
to be the multiplication by $\omega$ for $k=0,\cdots, 2n-2$.

The following two famous Lefschetz theorems give a filtration of
the Hodge groups and thus are extremely important in defining the
classifying space and the period map.

\begin{theorem} (Hard Lefschetz theorem) On a polarized algebraic variety
$(X,\omega)$ of dimension $n$,
\[
 L^k: H^{n-k}(X,\C)\rightarrow H^{n+k}(X,\C).
\]
is an isomorphism for every positive integer $k\leq n$.
\end{theorem}

The primitive cohomology $P^k(X,\C)$ is  then defined to be
the kernel of $L^{n-k+1}$ on $H^k(X,\C)$.

\begin{theorem}(Lefschetz Decomposition Theorem)
On a polarized algebraic variety $(X,\omega)$ of dimension $n$, we
have  the following decomposition:
\[
{\displaystyle H^n(X,\C)=\underset{k=0}{\overset{{[\frac
n2]}}{\oplus}} L^kP^{n-2k}(X,\C)}.
\]
\end{theorem}

Let  $H_Z=P^n(X,\C)\cap H^n(X,\mathbb Z)$ and  
$H^{p,q}=P^n(X,\C)\cap H^{p,q}(X)$ for $0\leq p,q\leq n$. Then we
have
\[
H_Z\otimes \C=\sum H^{p,q}, \quad H^{p,q}=\bar{H^{q,p}}
\]
for $p+q=n$. Set $H=H_Z\otimes \C$. We call $\{H^{p,q}\}$ the Hodge
decomposition of $H$.

\begin{rem}
We  define a filtration of $H_Z\otimes \C=H$ by
\[
0\subset F^n\subset F^{n-1}\subset \cdots F^1=H
\]
such that
\[
H^{p,q}=F^p\cap\bar{F}^q,\qquad F^p\oplus\bar{F^{n-p+1}}=H.
\]
The set $\{H^{p,q}\}$ and $\{F^p\}$ are equivalent in describing
the Hodge decomposition of $H$ (cf. \eqref{hdcomp}). We  will use both notations
interchangeably for the rest of this paper.
\end{rem}

Now suppose that $Q$ is the quadratic form on $H_Z$ 
induced by the cup product of the cohomology group
$H^n_{DR}(X,\C)$.  
$Q$ can be represented by
\begin{equation}\label{polar1}
Q(\phi, \psi)=(-1)^{n(n-1)/2}\int_X \phi\wedge \psi
\end{equation}
for $\phi, \psi\in H$.
$Q$ is a nondegenerat quadratic form, and is skew-symmetric if $n$
is odd  and is symmetric if $n$ is even. On $H$, the form $Q$ satisfies
the two  Hodge-Riemann relations on the space $H^{p,q}$ of primitive
harmonic $(p,q)$ forms: 
\begin{enumerate}\label{hr}
\item  $Q(H^{p,q},H^{p',q'})=0 \quad unless\quad p'=n-p,
q'=n-q$;
\item $(\sqrt{-1})^{p-q}\,Q(\phi,\bar\phi)>0$ {\it for any
nonzero  element} $\phi\in H^{p,q}$.
\end{enumerate}

\begin{definition}
A polarized Hodge structure
of weight $n$, denoted by
$\{H_Z,F^p,Q\}$, is given by a lattice $H_Z$, a 
filtration of $H=H_Z\otimes \C$
\[
0\subset F^n\subset F^{n-1}\subset\cdots\subset F^0\subset H,
\]
such that
\[
H=F^p\oplus\bar{F^{n-p+1}},
\]
together with a bilinear form
\[
Q: H_Z\otimes H_Z\rightarrow \Z,
\]
which is skew-symmetric if $n$ is odd and symmetric if $n$ is 
even such that it satisfies the two Hodge-Riemann
relations:

{\rm 3.}  $Q(F^p,F^{n-p+1})=0$ for $p=1,\cdots n$;

{\rm 4.}  $(\sqrt{-1})^{p-q}\,Q(\phi,\bar\phi)>0$ if $\phi\in 
H^{p,q}$ and $\phi\neq 0$,

\noindent where $H^{p,q}$ is defined by
\[
H^{p,q}=F^p\cap\bar{F^q}
\]
for $p+q=n$.
\end{definition}

\begin{definition}\label{def22}
The classifying space $D$ for the 
polarized Hodge structure is the set of all filtrations
\[
0\subset F^n\subset\cdots\subset F^1\subset H, \qquad
F^p\oplus\bar{F^{n-p+1}}=H,
\]
or the set of all the decompositions
\[
\sum H^{p,q}=H,\qquad H^{p,q}=\bar{H^{q,p}}
\]
on which $Q$ satisfies the two
Hodge-Riemann relations 1,2 or 3,4 above.
\end{definition}

Let
\begin{equation}\label{df23}
G_\R=\{\xi\in {\rm Hom}(H_\R,H_\R)|Q(\xi \phi,\xi\psi)=Q(\phi,\psi)\}.
\end{equation}
Then $D$ can also be written as the homogeneous space 
\begin{equation}\label{df24}
D=G/V,
\end{equation}
where $V$ is the compact subgroup of $G$ which leaves a fixed Hodge
decomposition $\{H^{p,q}\}$ invariant. Apparently, $G$ is a
semisimple real Lie group. 

Over the classifying space $D$ we have the holomorphic vector bundles
$\underline{F}^n,
\cdots, \underline{F}^1, \underline{H}$ whose fibers at each point are the
vector spaces $F^n,\cdots, F^1,H$, respectively. These bundles are called
Hodge bundles.

In Section~\ref{4}, we  identify the holomorphic tangent
bundle
$T^{1,0}(D)$  as a subbundle of $Hom(\underline
H,\underline H)$:
\[
T^{1,0}(D)\subset \oplus Hom(\underline F^p,\underline H/\underline
F^p)=
\underset{r> 0}{\oplus}
Hom( \underline H^{p,q},\underline H^{p-r,q+r}),
\]
such that the following compatible condition holds
\[
\begin{CD}
F^{p} @> >>F^{p-1}\\
@VVV        @VVV\\
H/F^p@ <<< H/F^{p-1}
\end{CD}.
\]

\begin{definition}  
A subbundle $T_h(D)$  is called the
horizontal distribution of $D$, if
\[
T_h(D)=\{\xi\in T^{1,0}(D)| \xi F^p\subset F^{p-1}, p=
1,\cdots, n\}.
\]
\end{definition}

For any point $x\in D$ such that $x$ is defined as subspaces 
$\{H^{p,q}\}$ of $H$, define the two vector spaces
\begin{align*}
& H^+=H^{n,0}+H^{n-2,2}+\cdots;\\
& H^-=H^{n-1,1}+H^{n-3,3}+\cdots.
\end{align*}
We fix a point $x_0\in D$. Suppose the corresponding vector spaces
are 
$\{H_0^{p,q}\}$ and $\{H_0^+, H_0^-\}$. Define $K$ to be the connected compact
subgroup of $G$ leaving $H_0^+$
invariant. We
give the basic properties of the classifying spaces in the following
three lemmas. 
The proofs are easy and are omitted.

\begin{lemma}
\label{121}
$K$ is the maximal compact subgroup of 
$G$ containing $V$. In particular, $V$ itself is a compact subgroup.
\end{lemma}
\qed

Define the Weil operator
\[
C: H^{p,q}\rightarrow H^{p,q},\quad C|_{H^{p,q}}=(\sqrt{-1})^{p-q}.
\]
Then we have
\[
C|_{H^+}=(\sqrt{-1})^n,\qquad C|_{H^-}=-(\sqrt{-1})^n.
\]
 
Let
\[
Q_1(x,{y})=Q(Cx,\bar{y}).
\]
Then we have

\begin{lemma}\label{lem122}
$Q_1$ is an Hermitian inner product.
\end{lemma} \qed
\begin{lemma}
Let
\[
D_1=\{H^{n,0}+H^{n-2,2}+\cdots|\{H^{p,q}\}\in D\}.
\]
Then the group $G$ acts on $D_1$ transitively with the stable subgroup 
$K$ at $H_0^+$, and  $D_1$ is a (Riemannian) symmetric space.
\end{lemma}
\qed

\begin{definition}\label{def24}
We call map $p$
\[
p: G/V\rightarrow G/K,\qquad \{H^{p,q}\}\mapsto H^{n,0}+H^{n-2,2}+\cdots
\]
the natural projection of the classifying space.
Using the notation of coset, $p(aV)=aK$ for any $a\in G$.
\end{definition}

With the above discussions, we can prove

\begin{prop}\label{pp21}
Suppose $T_v(D)$ is the distribution of the tangent vectors of the fibers 
of the canonical map
\[
p:\, D\rightarrow G/K,
\]
then
\[
T_v(D)\cap T_h(D)=\{ 0\}.
\]
\end{prop}

{\bf Proof:}
Let ${\frk g}$ be the Lie algebra of the Lie group $G$. Let ${\frk
g}={\frk f}+{\frk p}$ be the Cartan decomposition such that ${\frk f}$ is
the Lie algebra of $K$. Then
\[
T_v(D)=G\times_V{\frk v}_1
\]
where ${\frk f}={\frk v}+{\frk v}_1$ and ${\frk v}_1$ is the orthonormal 
complement of the Lie algebra ${\frk v}$ of $V$. On the other
hand, $T_h(D)\subset G\times_V{\frk p}$. So we have 
$T_v(D)\cap T_h(D)=\{ 0\}$.

 \qed

\begin{definition}
A horizontal slice ${\mathcal M}$ of $D$ is a complex integral submanifold
of the distribution $T_h(D)$.
\end{definition}

\begin{definition}
Let $U$ be an open neighborhood of the universal deformation 
space 
of $X$. Assume that $U$ is smooth. Then for each $X'$ near $X$, we
have an isomorphism $H^n(X',\C)=H^n(X,\C)$. Under this isomorphism, 
$\{H^{p,q}(X')\cap P^n(X',\C)\}_{p+q=n}$ can be considered as a
point of
$D$. The map 
\[
U\rightarrow D,\qquad X'\mapsto \{H^{p,q}(X')\cap
P^n(X',\C)\}_{p+q=n}
\]
is called the period map. If $\Gamma_1\rightarrow\Gamma$ is a
homomorphism between two discrete groups and the period map is
equivariant with respect to the two groups, then we also call
induced map
\[
\Gamma_1\backslash U\rightarrow\Gamma\backslash D
\]
a period map.
\end{definition}

The most important property of the period map is the following~\cite{Gr}:

\begin{theorem} (Griffiths)
The period map $p: U\rightarrow D$ is holomorphic. Furthermore, it is  an
immersion and  
$p(U)$ is 
a  horizontal slice of the classifying space.
\end{theorem}

From the above theorem and Proposition~\ref{pp21} in this section, we can
prove:

\begin{cor}\label{cor31}
With the notations as above, the map
\[
p: U\subset D\rightarrow D_1=G/K
\]
is a (real) immersion.
\end{cor}

\begin{definition}\label{hodgemetric}
Using the above notations, let $h$ be the invariant \ka metric on
$D_1$. The Hodge metric is defined as  the (Riemannian) metric $p^*h$ on the
horizontal slice $U$.
\end{definition}

\begin{rem}
In ~\cite{Lu3}, the first author proved that the Hodge metric of $U$ is \ka. 
\end{rem}

Now we introduce the Nilpotent Orbit theorem of
Schmid~\cite{Schmid}. Let
$f: {\mathcal X}\rightarrow S$ be a family of compact \ka manifolds.
In order to study the degeneration of the variation of the Hodge
structure, we let $S=\Delta^{*l}\times \Delta^{m-l}$,
where $l\geq 1, m\geq l$, and $\Delta$, $\Delta^*$ are the unit
disk and the punctured unit disk in the complex plane, respectively.
Consider the  period map
\[
\Phi:\Delta^{*l}\times \Delta^{m-l}\rightarrow \Gamma\backslash D.
\]
 By going to the universal covering
$U^l\times \Delta^{m-l}$, one can lift $\Phi$ to a mapping
\[
\tilde\Phi: U^l\times \Delta^{m-l}\rightarrow D,
\]
where $U$ is the upper half plane.
Corresponding to each of the first $l$ variables, we choose a
monodromy transformation $T_i\in \Gamma$, where $\Gamma$ is the
monodromy group,  so that
\[
\tilde\Phi(z_1,\cdots, z_i+1,\cdots z_l,w_{l+1},\cdots,w_m)
=T_i\circ\tilde\Phi(z_1,\cdots,z_l,w_{l+1},\cdots,w_m),
\]
holds identically in all variables. $T_i$'s commute
with each other. We know that all the eigenvalues of $T_i$ are roots 
of unity. Let $T_i=T_{i,s}T_{i,u}$ be the Jordan decomposition where 
$T_{i,s}$ is semisimple and $T_{i,u}$ is unipotent. 
We also assume that $T_{i,s}^{s_{i}}=I$ for some
positive integer $s_i$ so that we can define
$N_i=\frac{1}{s_i}\log T_i^{s_i}=\sum_{k \geq 1}
(-1)^{k+1}\frac{1}{k}(T_i^{s_i}-I)^k$. All $N_i$ are  commutative. 

Let $z=(z_1,\cdots,z_l)$, $sz=(s_1z_1,\cdots,s_lz_l)$ and
$w=(w_{l+1},\cdots,w_m)$. The map
\[
\tilde\Psi(z,w)=\exp(-\sum_{i=1}^ls_iz_iN_i)\circ\tilde\Phi(sz,w)
\]
remains invariant under the translation $z_i\mapsto z_i+1, 1\leq
i\leq l$. It follows that $\tilde\Psi$ drops to a mapping
\[
\Psi:\Delta^{*l}\times\Delta^{m-l}\rightarrow \check{D}.
\]

\begin{theorem}[\bf Nilpotent Orbit Theorem~\cite{Schmid}]\label{tnot}
The map
$\Psi$ extends holomorphically to $\Delta^m$. For $w\in
\Delta^{m-l}$, the point
\[
a(w)=\Psi(0,w)\in\check D
\]
is left fixed by $T_{i,s}, 1\leq i\leq l$. 
For any given number $\eta$ with $0<\eta<1$, there exist constants
$\alpha,\beta\geq 0$, such that under the restrictions
\[
{\rm Im} z_i\geq\alpha, 1\leq i\leq l\quad{\rm and}\quad
|w_j|\leq\eta, l+1\leq j\leq m,
\]
the point $exp(\sum_{i=1}^lz_iN_i)\circ a(w)$ lies in $D$ and
satisfies the inequality
\[
d(exp(\sum_{i=1}^lz_iN_i)\circ a(w), \tilde\Phi(z,w))
\leq
(\Pi_{i=1}^l{\rm Im} z_i)^\beta\sum_{i=1}^l exp(-2\pi {\rm s_{i}^{-1}Im }z_i)
\]
here $d$ is the $G_R$ invariant Riemannian distance function on
$D$. Finally, the mapping
\[
(z,w)\mapsto exp(\sum_{i=1}^lz_iN_i)\circ a(w)
\]
is horizontal. 
\end{theorem} 

Now we assume that the generic fiber $X$ of the map $f:{\mathcal
X}\rightarrow S$ is a polarized Calabi-Yau manifold.  For the the
sake of simplicity, we assume that $X$ is compact, simply connected
and algebraic with
$c_1(X)=0$. By a theorem of Tian~\cite{T1}
\footnote{Tian's proof is more general since one merely assumes the
$\pa\bar\pa$-lemma hold for $X$. That is equivalent to assume that 
the Hodge-de Rham spectral sequence for $X$ degenerates at the $E_1$ term.
See the survey paper of Friedman~\cite{Fd} for details. }
 (See also
Todorov~\cite{To}), the universal deformation space of
$X$ is smooth. Since there are no nonzero holomorphic vector fields
on a Calabi-Yau manifold, the 
moduli space of polarized Calabi-Yau manifolds is an orbifold.
The following important theorem of Viehweg gives the
compactification of  the moduli space:

\begin{theorem} ({\bf Viehweg~\cite[page 21, Theorem 1.13]{V2}})
Let ${\mathcal M}$ be the moduli space of polarized Calabi-Yau
manifolds and the line bundle $\underline{F}^n$ is the Hodge bundle
defined right after Definition~\ref{def22}. Then ${\mathcal M}$ is
quasi-projective and the line bundle $\underline{F}^n$ extends to an
ample line bundle over
$\bar{\mathcal M}$, the compactification of ${\mathcal M}$.
\end{theorem}

With the classical Hironaka theorem~\cite{Hi}, we have the following
\begin{cor}
Let $\bar{\mathcal M}$ be the compactification of the moduli space 
${\mathcal M}$ in the above sense, then after a smooth resolution,
 one can assume that 
$\bar{\mathcal M}\backslash{\mathcal M}$ is a divisor of normal
crossing. In other word, let $x_0\in \bar{\mathcal
M}\backslash{\mathcal M}$, then in a neighborhood of $x_0$, we can 
write ${\mathcal M}$ as
\[
\Delta^{*l}\times (\Delta)^{m-l},
\]
where $m$ is the complex dimension of $\bar{\mathcal
M}$. 
\end{cor}

\section{Curvature of Weil-Petersson metrics}\label{S3}
For the rest of this paper, we assume that ${\mathcal M}$ is the moduli space of
polarized Calabi-Yau manifolds of dimension $n>2$~\footnote{For $K3$
surfaces, the Weil-Petersson metric is half of the Hodge metric. Thus we omit
this case.}.
 
\begin{rem}\label{setconvention}
The following notations and conventions will be used through out the rest of
this paper. 

Form a \ka manifold $M$ with metric $g_{i\bar j}$, the curvature tensor is given
by 
\[
R_{i\bar j k\bar l}=\frac{\partial^2 g_{i\bar j}}{\partial z_k\partial\bar{z_l}}
-g^{p\bar q}\frac{\partial g_{i\bar q}}{\partial z_k}
\frac{\partial g_{p\bar j}}{\partial\bar{z_l}}.
\]
Using this convention, the Ricci curvature is 
\[
R_{i\bar j}=-g^{k\bar l}R_{i\bar j k\bar l}.
\]
Furthermore, the Christoffel symbol of this metric is given by
\[
\Gamma_{ij}^k=g^{k\bar q}\frac{\pa g_{i\bar q}}{\pa z_j}.
\]
We also have
\[
\frac{\pa \Gamma_{ij}^k}{\pa \bar{z_j}}=g^{k\bar q}R_{i\bar lj\bar q}.
\]
So the holomorphic bisectional curvature of this metric is nonpositive means 
$R_{i\bar i j\bar j}\geq 0$ for all $i,j$.
\end{rem}

Now let $\Omega$ be a nonzero local holomorphic section of the Hodge 
bundle 
$\underline{F}^{n}$. In this section, we assume $1 \leq i,j \leq m$ unless 
otherwise stated, where $m$ is the dimension of the moduli space ${\mathcal M}$. We
set  ~\footnote{In fact, we use the notation $(\xi,\eta)=(\sqrt{-1})^nQ(\xi,\eta)$
in the rest of this paper where $\xi,\eta$ are $n$-forms. The bilinear form
$(\,,\,)$ is not necessary positive definite.}
\begin{equation}\label{scalar}
(\Omega,\bar{\Omega})=(\sqrt{-1})^{n}Q(\Omega,\bar{\Omega}).
\end{equation}
By the Hodge-Riemann relations we know that $(\Omega,\bar{\Omega})>0$. 
In local coordinates, the Weil-Petersson metric is given by 
\begin{eqnarray}\label{defwp}
g_{i\bar{j}}=-\partial_{i}\bar\partial_{{j}} \log (\Omega, \bar{\Omega})
=-\frac{(\partial_{i}\Omega, \bar{\partial_{j}\Omega})}{(\Omega, \bar{\Omega})}
+\frac{(\partial_{i}\Omega, \bar{\Omega})(\Omega,
\bar{\partial_{j}\Omega})}{(\Omega, \bar{\Omega})^{2}},
\end{eqnarray}
where $\pa_i,\bar\pa_j$ are the operators
$\frac{\pa}{\pa z_i},\frac{\pa}{\pa \bar z_j}$, respectively.
From~\cite{T1}, we know that the definition is the same as the Weil-Petersson metric
defined in the classical way. In this section, we compute the curvature
of the Weil-Petersson metric. We begin with defining
\begin{eqnarray}\label{*1}
K_{i}=-\partial_{i} \log (\Omega, \bar{\Omega})=-\frac{(\partial_{i}\Omega, 
\bar{\Omega})}{(\Omega, \bar{\Omega})},
\end{eqnarray}
and
\begin{eqnarray}\label{*2}
D_{i}\Omega=\partial_{i}\Omega+K_{i}\Omega
\end{eqnarray}
for $1 \leq i \leq m$. Then $g_{i\bar{j}}=\bar\partial_{{j}}K_{i}$.
\begin{lemma}\label{h31}
Under the notions as above, the following properties hold:

(1) $(D_{i}\Omega,\bar{\Omega})=0$; 

(2) $\bar\partial_{{j}} D_{i}\Omega=g_{i\bar{j}} \Omega$; 

(3) $g_{i\bar{j}}=-\frac{(D_{i}\Omega,\bar{D_{j}\Omega})}
{(\Omega, \bar{\Omega})}$, 

\noindent where $1 \leq i,j \leq m$.
\end{lemma}

{\bf Proof.}
By \eqref{*1} and \eqref{*2} we have 
\[
(D_{i}\Omega,\bar{\Omega})=(\partial_{i}\Omega+K_{i}\Omega,\bar{\Omega})
=(\partial_{i}\Omega,\bar{\Omega})-\frac{(\partial_{i}\Omega, 
\bar{\Omega})}{(\Omega, \bar{\Omega})}(\Omega, \bar{\Omega})=0,
\]
which proves (1). (2) follows from 
\[
\bar\partial_{{j}} D_{i}\Omega=\bar\partial_{{j}}\partial_{i}\Omega
+(\bar\partial_{{j}}K_{i})\Omega=g_{i\bar{j}}\Omega.
\]
Combining the above two equations with \eqref{defwp} we have 
\begin{eqnarray*}
g_{i\bar{j}} = \frac{(\bar\partial_{{j}} D_{i}\Omega, \bar{\Omega})}
{(\Omega,\bar{\Omega})}=-\frac{(D_{i}\Omega,\bar{\partial_{j}\Omega})}
{(\Omega, \bar{\Omega})}=-\frac{(D_{i}\Omega,\bar{D_{j}\Omega})}
{(\Omega, \bar{\Omega})}+\frac{(D_{i}\Omega,\bar{K_{j}\Omega})}
{(\Omega, \bar{\Omega})}
= -\frac{(D_{i}\Omega,\bar{D_{j}\Omega})}
{(\Omega, \bar{\Omega})}.
\end{eqnarray*}
This finishes the proof.

\qed

From the above lemma, we see that $D_{i}\Omega$ is the projection of 
$\partial_{i}\Omega$ into $H^{n-1,1}$
with respect to the quadratic form $(\,,\,)$.
 Now we consider the projection 
of $\partial_{j}D_{i}\Omega$ into $H^{n-2,2}$. In the following we will use 
$\Gamma_{ij}^k$ to denote the Christoffel symbol of the Weil-Petersson metric. 
Let 
\begin{eqnarray}\label{*3}
D_{j}D_{i}\Omega=\partial_{j}D_{i}\Omega - \sum_{k} \Gamma^{k}_{ij}D_{k}\Omega
+K_{j}D_{i}\Omega.
\end{eqnarray}
\begin{lemma}\label{h22}
Using the same notations as above, for any $1 \leq i,j,l \leq m$, we have 

(1) $(D_{j}D_{i}\Omega,\bar{\Omega})=0$;

(2) $(D_{j}D_{i}\Omega,\bar{D_{l}\Omega})=0$;

(3) $D_{j}D_{i}\Omega=D_{i}D_{j}\Omega$.

\end{lemma}

{\bf Proof.}
A straightforward computation gives 
\begin{eqnarray*}
(D_{j}D_{i}\Omega,\bar{\Omega})=(\partial_{j}D_{i}\Omega,\bar{\Omega})
-\sum_{k} \Gamma^{k}_{ij}(D_{k}\Omega,\bar{\Omega})+
K_{j}(D_{i}\Omega,\bar{\Omega})
=\partial_{j}(D_{i}\Omega,\bar{\Omega})=0,
\end{eqnarray*}
where
in the last equality, we used (1) of Lemma~\ref{h31}.
This proves (1). Using Lemma~\ref{h31}, we have 
\begin{eqnarray*}
& & (D_{j}D_{i}\Omega,\bar{D_{l}\Omega})
=(\partial_{j}D_{i}\Omega,\bar{D_{l}\Omega})
-\sum_{k} \Gamma^{k}_{ij}(D_{k}\Omega,\bar{D_{l}\Omega})
+K_{j}(D_{i}\Omega,\bar{D_{l}\Omega})\\
&=& \partial_{j}(-g_{i\bar{l}}(\Omega,\bar{\Omega}))
-(D_{i}\Omega,\bar{\bar\partial_{{j}} D_{l}\Omega}) +\sum_{k}
\Gamma^{k}_{ij}g_{k\bar{l}}(\Omega,\bar{\Omega})
-K_{j}g_{i\bar{l}}(\Omega,\bar{\Omega})\\
&=& -\partial_{j}g_{i\bar{l}}(\Omega,\bar{\Omega})
-g_{i\bar{l}}(\partial_{j}\Omega,\bar{\Omega})
-(D_{i}\Omega,g_{j\bar{l}}\bar{\Omega}) +\partial_{j}g_{i\bar{l}}
(\Omega,\bar{\Omega})
-g_{i\bar{l}}(K_{j}\Omega,\bar{\Omega})\\
&= &-g_{i\bar{l}}(D_{j}\Omega,\bar{\Omega})-g_{j\bar{l}}
(D_{i}\Omega,\bar{\Omega})
=0.
\end{eqnarray*}
This proves (2). 
To prove (3), we see that 
\begin{eqnarray*}
D_{j}D_{i}\Omega &=&\partial_{j}D_{i}\Omega - \sum_{k} \Gamma^{k}_{ij}D_{k}\Omega
+K_{j}D_{i}\Omega \\
&=&\partial_{j}\partial_{i}\Omega
+K_{i}\partial_{j}\Omega 
- \sum_{k} \Gamma^{k}_{ij}D_{k}\Omega +K_{j}\partial_{i}\Omega +K_{j}K_{i}\Omega 
-\frac{(\partial_{j}\partial_{i}\Omega,\bar{\Omega})}
{(\Omega,\bar{\Omega})}\Omega\\
& &+\frac{(\partial_{i}\Omega,\bar{\Omega})(\partial_{j}\Omega,\bar{\Omega})}
{(\Omega,\bar{\Omega})^{2}}\Omega. \hspace{2.2in}
\end{eqnarray*}
Thus (3) follows from the fact that the above formula is symmetric with
respect to $i$ and $j$. 

\qed

Let
$R_{i\bar{j}k\bar{l}}$ be the  curvature tensor of $g_{i\bar{j}}$. Then we have
the following~\cite{Wang1}:
\begin{theorem}\label{Str}
Let $(g_{i\bar j})_{m\times m}$ be the Weil-Petersson metric and let 
$D_{j}D_{i}\Omega$ be defined as in \eqref{*3}. Then the Weil-Petersson metric 
is \ka \cite{T1}, and the curvature tensor is  
\begin{eqnarray}\label{strominger}
R_{i\bar{j}k\bar{l}}=g_{i\bar{j}}g_{k\bar{l}}+g_{i\bar{l}}g_{k\bar{j}}
-\frac{(D_{k}D_{i}\Omega,\bar{D_{l}D_{j}\Omega})}{(\Omega,\bar{\Omega})}
\end{eqnarray}
for $1 \leq i,j,k,l \leq m$. 
\end{theorem}

{\bf Proof.}
By definition 
\begin{eqnarray}\label{curv}
R_{i\bar{j}k\bar{l}}=\frac{\partial^{2}g_{i\bar{j}}}
{\partial z_{k} \partial \bar z_{{l}}}-g^{p\bar{q}}
\frac{\partial g_{i\bar{q}}}{\partial z_{k}}
\frac{\partial g_{p\bar{j}}}{\partial \bar z_{{l}}}. 
\end{eqnarray}
From (3) of Lemma~\ref{h31},  we know 
\[
\frac{\partial g_{i\bar{j}}}{\partial z_{k}}=
-\frac{(\partial_{k}D_{i}\Omega,\bar{D_{j}\Omega})}{(\Omega,\bar{\Omega})}
-\frac{(D_{i}\Omega,\partial_{k}\bar{D_{j}\Omega})}{(\Omega,\bar{\Omega})}
+\frac{(D_{i}\Omega,\bar{D_{j}\Omega})}{(\Omega,\bar{\Omega})^{2}}
(\partial_{k}\Omega, \bar{\Omega}).
\]
By (1) of Lemma~\ref{h31}, 
\[
\frac{(D_{i}\Omega,\partial_{k}\bar{D_{j}\Omega})}{(\Omega,\bar{\Omega})}
=\frac{(D_{i}\Omega,\bar{\bar\partial_{{k}}D_{j}\Omega)}}{(\Omega,\bar{\Omega})}
=\frac{(D_{i}\Omega,g_{k\bar{j}}\bar{\Omega})}{(\Omega,\bar{\Omega})}=0.
\]
Using the definition of $K_i$,
we have 
\begin{eqnarray*}
\frac{\partial g_{i\bar{j}}}{\partial z_{k}}
=-\frac{(\partial_{k}D_{i}\Omega+K_{k}D_{i}\Omega,\bar{D_{j}\Omega})}
{(\Omega,\bar{\Omega})}.
\end{eqnarray*}
Let $A_{ij}=\partial_{i}D_{j}\Omega
+K_{i}D_{j}\Omega=D_{i}D_{j}\Omega+
\Gamma^{k}_{ij}D_{k}\Omega$.
Then
\begin{eqnarray}\label{1st}
\frac{\partial g_{i\bar{j}}}{\partial z_{k}}=-\frac{(A_{ki},\bar{D_{j}\Omega})}
{(\Omega,\bar{\Omega})}.
\end{eqnarray}
Similarly, we have 
\begin{eqnarray}\label{1stbar}
\frac{\partial g_{i\bar{j}}}{\partial \bar z_{{l}}}
=-\frac{(D_{i}\Omega,\bar{A_{lj}})}{(\Omega,\bar{\Omega})}
.\end{eqnarray}
From \eqref{1st} we have 
\begin{eqnarray}\label{2nd}
\frac{\partial^{2}g_{i\bar{j}}}{\partial z_{k} \partial \bar z_{{l}}}
=-\frac{(\bar\partial_{{l}}A_{ki},\bar{D_{j}\Omega})}{(\Omega,\bar{\Omega})}
-\frac{(A_{ki},\bar{\partial_{l}D_{j}\Omega})}{(\Omega,\bar{\Omega})}
+\frac{(A_{ki},\bar{D_{j}\Omega})}{(\Omega,\bar{\Omega})^{2}}
(\Omega,\bar{\partial_{l}\Omega}). 
\end{eqnarray}
We also have 
\begin{eqnarray*}
\bar\partial_{{l}}A_{ki}
&=&\bar\partial_{{l}}(\partial_{k}D_{i}\Omega+K_{k}D_{i}\Omega)
=\partial_{k}(\bar\partial_{{l}}D_{i}\Omega)+(\bar\partial_{{l}}K_{k})D_{i}\Omega
+K_{k}\bar\partial_{{l}}D_{i}\Omega\\
&=&\partial_{k}(g_{i\bar{l}}\Omega)+g_{k\bar{l}}D_{i}\Omega+
K_{k}g_{i\bar{l}}\Omega
=(\partial_{k}g_{i\bar{l}})\Omega
+g_{i\bar{l}}D_{k}\Omega+g_{k\bar{l}}D_{i}\Omega,
\end{eqnarray*}
and
\begin{eqnarray*}
\frac{(A_{ki},\bar{D_{j}\Omega})}{(\Omega,\bar{\Omega})^{2}}
(\Omega,\bar{\partial_{l}\Omega})
=-\frac{(A_{ki},\bar{K_{l}D_{j}\Omega})}{(\Omega,\bar{\Omega})}.
\end{eqnarray*}
Thus from~\eqref{2nd}, we have
\begin{align}\label{2nd2}
\begin{split}
\frac{\partial^{2}g_{i\bar{j}}}{\partial z_{k} \partial \bar z_{{l}}}
=&-\frac{((\partial_{k}g_{i\bar{l}})\Omega
+g_{i\bar{l}}D_{k}\Omega+g_{k\bar{l}}D_{i}\Omega,
\bar{D_{j}\Omega})}{(\Omega,\bar{\Omega})}
-\frac{(A_{ki},\bar{\partial_{l}D_{j}\Omega})}{(\Omega,\bar{\Omega})} 
-\frac{(A_{ki},\bar{K_{l}D_{j}\Omega})}{(\Omega,\bar{\Omega})}\\ 
=& g_{i\bar{l}}g_{k\bar{j}}+g_{i\bar{j}}g_{k\bar{l}}-
\frac{(A_{ki},\bar{A_{lj}})}{(\Omega,\bar{\Omega})},
\end{split}
\end{align}
by using Lemma~\ref{h31}.
Combining \eqref{1st}, \eqref{1stbar} and \eqref{2nd2}, and using
Lemma~\ref{h22},  we have 
\begin{eqnarray*}
R_{i\bar{j}k\bar{l}} &=& g_{i\bar{j}}g_{k\bar{l}}+g_{i\bar{l}}g_{k\bar{j}}-
\frac{(A_{ki},\bar{A_{lj}})}{(\Omega,\bar{\Omega})}+g^{p\bar{q}}
\frac{\pa g_{i\bar q}}{\pa z_k}
\frac{(D_{p}\Omega,\bar{A_{lj}})}{(\Omega,\bar{\Omega})} \\
&=& g_{i\bar{j}}g_{k\bar{l}}+g_{i\bar{l}}g_{k\bar{j}}-
\frac{(A_{ik}-\Gamma^{p}_{ik}D_{p}\Omega, \bar{A_{lj}})}{
(\Omega,\bar\Omega)}\\
&=& g_{i\bar{j}}g_{k\bar{l}}+g_{i\bar{l}}g_{k\bar{j}}-
\frac{(D_{k}D_{i}\Omega,
\bar{D_{l}D_{j}\Omega})}{(\Omega,\bar{\Omega})}.
\end{eqnarray*}
This finishes the proof.
\qed

\begin{rem}
For the moduli space of Calabi-Yau threefolds, Strominger \cite{S} proved that 
the curvature tensor is 
\begin{eqnarray*}
R_{i\bar j k \bar l}=g_{i\bar j}g_{k \bar l}+g_{i\bar l}g_{k \bar j}
-\sum_{p,q}\frac{1}{(\Omega,\bar\Omega)^2}g^{p\bar q} F_{ikp}\bar{F_{jlq}},
\end{eqnarray*}
where $F_{ikp}$ is the Yukawa coupling. In the case of Calabi-Yau threefolds, 
$D_{i}D_{j}\Omega \in H^{1,2}$. 
In fact, $D_iD_j\Omega$ is the orthogonal projection of
$\pa_i\pa_j\Omega$ to $H^{1,2}$.
Thus
\begin{eqnarray*}
(D_{i}D_{k}\Omega, \bar{D_{j}D_{l}\Omega})
=-\frac{(D_{i}D_{k}\Omega, D_{p}\Omega)
(\bar{D_{q}\Omega}, \bar{D_{j}D_{l}\Omega})}{(\Omega,\bar\Omega)}g^{p\bar q}.
\end{eqnarray*}
It is easy to see that $(D_{i}D_{k}\Omega, D_{p}\Omega)=-F_{ikp}$. Thus our 
theorem is the same as Strominger's in the case of Calabi-Yau threefolds. 
\end{rem}

\begin{rem}
Theorem~\ref{Str} was proved in
\cite{Sc1} using the method of
\cite{Siu2}  which is different from ours. 
In his paper~\cite{To}, Todorov introduced the geodesic coordinates from
which it is much easier to get the curvature formula. The current proof
was by
Wang~\cite{Wang1} which is purely  Hodge theoretic. Such a proof
can be generalized to  general horizontal slice.  
\end{rem}

\section{Partial Hodge Metrics}\label{3half}
We use the Ricci curvature of the Weil-Petersson metric to construct a
new  metric. Let $\omega_{WP}$ be the \ka form of the Weil-Petersson metric
and let
$\mu >m+1$ be a real number. Let 
\begin{eqnarray}\label{newmetric}
\omega_{\mu}=\mu \omega_{WP}+Ric(\omega_{WP}).
\end{eqnarray}
By Theorem~\ref{Str}, we know that $\omega_\mu$ is a \ka metric.
We notice here that when the dimension of the Calabi-Yau 
manifolds is $3$ or $4$, by choosing suitable $\mu$, the metric $\omega_{\mu}$ 
 coincides with the Hodge metric(cf. \S~\ref{4}). 
For this reason we call $\omega_\mu$ the ``partial Hodge metric''. It
is a metric between the Weil-Petersson metric and the Hodge metric.

In this section, unless otherwise stated, the subscripts
$i,j,k,l,p,q,\alpha,\cdots$ will be ranging from $1$ to $m$.
Define a tensor
\begin{eqnarray}\label{*4}
T_{k\alpha i}=\partial_{k}D_{\alpha}D_{i}\Omega+K_{k}D_{\alpha}D_{i}\Omega 
-\sum_{p}\Gamma_{\alpha k}^{p}D_{p}D_{i}\Omega
-\sum_p\Gamma_{ik}^{p}D_{\alpha}D_p\Omega,
\end{eqnarray}
where $\Gamma_{\alpha k}^{p}$ is the Christoffel symbol of the 
Weil-Petersson metric and $\Omega$ is a nonzero local holomorphic section of
$\underline{F}^n$ as in the previous section. We  use
$g_{i\bar j}$  and $h_{i\bar j}$ to denote  the metric matrices of the
Weil-Petersson metric and the metric
$\omega_{\mu}$ (for some chosen $\mu$) in local coordinates
$(z_{1},\cdot\cdot\cdot,z_{m})$ respectively, and use 
$R_{i\bar j k\bar l}$ and 
$\widetilde{R}_{i\bar j k\bar l}$ to denote their curvature tensors
respectively. We also use $R_{i\bar j}$ to denote the Ricci tensor of 
the Weil-Petersson metric. 

Let $D_{k}D_{\alpha}D_{i}\Omega$ be the projection of 
$T_{k\alpha i}$ into $H^{n-3,3}$ with respect to the quadratic form
$(\,,\,)$ in ~\eqref{scalar}. 
Let $E_{k\alpha i}=T_{k\alpha i}-D_{k}D_{\alpha}D_{i}\Omega$. 
Then we have the following

\begin{lemma}\label{lem91}
 Using the same notations as above, we have
\[
T_{k\alpha i}=E_{k\alpha i}+D_kD_\alpha D_i\Omega\in H^{n-2,2}\oplus
H^{n-3,3},
\]
where $T_{k\alpha i}$ is defined in ~\eqref{*4}.
\end{lemma}

{\bf Proof.} By definition of $T_{k\alpha i}$ and the Griffiths' transversality, 
\[
T_{k\alpha i}\in H^{n,0}\oplus H^{n-1,1}\oplus H^{n-2,2}\oplus H^{n-3,3}.
\]
Using Lemma~\ref{h22},  we have
\[
(T_{k\alpha i}, \bar\Omega)=(\pa_kD_\alpha D_i\Omega, \bar\Omega)=0.
\]
So there is no $H^{n,0}$ components in $T_{k\alpha i}$. On the other hand,
$H^{n-1,1}$ is spanned by $D_i\Omega$. Using Lemma~\ref{h22} again, we have
\[
(T_{k\alpha i}, \bar{D_j\Omega})=0.
\]
Thus $T_{k\alpha i}$ has no $H^{n-1,1}$ component and
this completes the proof.

\qed

Define the curvature like tensor $F$ by 
\begin{eqnarray}\label{yucawa}
F_{i\bar j k\bar l}=\frac{(D_{k}D_{i}\Omega,\bar{D_{l}D_{j}\Omega})}
{(\Omega,\bar\Omega)}.
\end{eqnarray}
Using Lemma \ref{h22} and the Hodge-Riemann relations we know that the 
tensor $F$ has all symmetries that a curvature tensor has. 

The Strominger formula (Theorem~\ref{Str}) can be written as 
\begin{eqnarray}\label{stro1}
R_{i\bar j k\bar l}=g_{i\bar j}g_{k\bar l}+
g_{i\bar l}g_{k\bar j}-F_{i\bar j k\bar l}.
\end{eqnarray}

The curvature tensor of the partial Hodge metric is
\begin{theorem}\label{curvmu}
The metric $\omega_{\mu}$ is \ka and the curvature tensor of 
$\omega_{\mu}$ is 
\begin{align}\label{gencurv}
\begin{split}
\widetilde{R}_{i\bar j k\bar l} =& (\mu -m-1)(g_{i\bar j}g_{k\bar l}+
g_{i\bar l}g_{k\bar j})-(\mu -m)F_{i\bar j k\bar l}
+\sum_{\alpha\beta pq}F_{i\bar q \alpha\bar l}F_{p\bar j k\bar\beta}
g^{\alpha\bar\beta}g^{p\bar q}\\
&+\sum_{\alpha\beta pq}F_{\alpha\bar q k\bar l}F_{i\bar j p\bar\beta}
g^{\alpha\bar\beta}g^{p\bar q}
+\sum_{\alpha\beta}\frac{(D_{k}D_{\alpha}D_{i}\Omega,
\bar{D_{l}D_{\beta}D_{j}\Omega})}{(\Omega,\bar\Omega)}g^{\alpha\bar\beta}\\
&+\sum_{\alpha\beta}\frac{(E_{k\alpha i},
\bar{E_{l\beta j}})}{(\Omega,\bar\Omega)}g^{\alpha\bar\beta}
-\sum_{\alpha\beta\gamma\tau st}h^{s\bar t}\frac{(E_{k\alpha
i},\bar{D_{\beta}D_{t}\Omega})}
{(\Omega,\bar\Omega)}\frac{(D_{\gamma}D_{s}\Omega,
\bar{E_{l\tau j}})}{(\Omega,\bar\Omega)}g^{\alpha\bar\beta}
g^{\gamma\bar\tau}.
\end{split}
\end{align}
\end{theorem}
\noindent
We will leave the proof of this theorem in the appendix due to its length.

The main theorem of this section is that, for the moduli space 
of Calabi-Yau fourfolds, with a suitable choice of $\mu$, the partial
Hodge metric has the following property:
 
\begin{theorem}\label{fourdim}
Let $n=4$ and let $\omega_{PH}=(m+2)\omega_{WP}+Ric(\omega_{WP})$, then\\
1. $\omega_{PH}$ is a \ka metric.\\
2. The Ricci curvature and the holomorphic sectional curvature of 
$\omega_{PH}$ are bounded above by\\ 
\hspace*{0.1in}  the negative constant $-\frac{1}{m+4}$. \\
3. The holomorphic bisectional curvature of $\omega_{PH}$ is nonpositive.
\end{theorem}

{\bf Proof.}
By Theorem \ref{curvmu} we know that $\omega_{PH}$ is \ka. 
From \eqref{hodge} we know that 
\begin{eqnarray}\label{fourhodge}
h_{i\bar{j}}=
g_{i\bar{j}}+F_{i\bar j\alpha\bar\beta}g^{\alpha\bar\beta}.
\end{eqnarray} 
Fix  a point $x_0$ in the moduli space. Let 
$z_{1},\cdot \cdot \cdot ,z_{m}$ be the local holomorphic 
normal coordinate at
$x_0$ with respect to the Weil-Petersson metric. Then at the point $x_0$, we
have 
\begin{align}\label{ass3}
\begin{split}
&g_{\alpha\bar\beta}=\delta_{\alpha\beta},\\
&\Gamma^{\gamma}_{\alpha\beta}=\frac{\partial g_{\alpha\bar\beta}}{\partial
z_{\gamma}} =\frac{\partial g^{\alpha\bar\beta}}{\partial
z_{\gamma}}=\frac{\partial g_{\alpha\bar\beta}}{\partial
\bar z_{\gamma}}=\frac{\partial g^{\alpha\bar\beta}}{\partial
\bar z_{\gamma}}=0.
\end{split}
\end{align}
Replacing $\Omega$ by $\widetilde{\Omega}=f\Omega$,  where  $f$ is a local
holomorphic function defined by
\begin{eqnarray*}
f(z)=(\Omega,\bar{\Omega})^{-\frac{1}{2}}(x_0)-\sum_{i}
\frac{(\partial_{i}\Omega,\bar\Omega)(x_0)}{((\Omega,\bar\Omega)(x_0))
^{\frac 32}}z_{i}\ ,
\end{eqnarray*} 
we have , at  the point $x_0$, 
\begin{eqnarray}\label{ass1}
(\partial_{k}\tilde\Omega,\bar{\tilde\Omega})=(\tilde\Omega,
\bar{\partial_{k}\tilde\Omega})=0 
\end{eqnarray}
for each $k=1,2,\cdot \cdot \cdot ,n$ 
and 
\begin{eqnarray}\label{ass2}
(\tilde\Omega,\bar{\tilde\Omega})=1 .
\end{eqnarray}
By abusing of notations, we use $\Omega$ to replace $\tilde\Omega$ 
for the rest of this section.

We set $i=j$ and 
$k=l$.  Based on the above notations, from Theorem \ref{curvmu}  
 we have
\begin{align}\label{h2nd3}
\begin{split} 
&\widetilde{R}_{i\bar{i}k\bar{k}}
= 1+\delta_{ik}-2(D_{k}D_{i}\Omega,\bar{D_{k}D_{i}\Omega})
+\sum_{\alpha,\gamma}(D_{\alpha}D_{i}\Omega,\bar{D_{k}D_{\gamma}\Omega})
(D_{k}D_{\gamma}\Omega, \bar{D_{\alpha}D_{i}\Omega}) \\
&+\sum_{\alpha,\beta}(D_{\alpha}D_{i}\Omega,\bar{D_{\beta}D_{i}\Omega})
(D_{\beta}D_{k}\Omega,\bar{D_{\alpha}D_{k}\Omega}) 
+\sum_{\alpha} (D_k D_\alpha D_i\Omega, \bar{D_k D_\alpha D_i\Omega}) \\
&+\bigg ( \sum_{\alpha}(E_{k\alpha i}, \bar{E_{k\alpha i}})  
-h^{p\bar q}
\sum_{\alpha}(E_{k\alpha i},\bar{D_{\alpha}D_{q}\Omega})
\sum_{\beta}(D_{\beta}D_{p}\Omega,\bar{E_{k\beta i}})
\bigg ). 
\end{split}
\end{align}

Fix the indices $k$ and $i$. Let $U_{\alpha}=h^{p\bar q}
\sum_{\beta}(E_{k\beta i},\bar{D_{\beta}D_{q}\Omega})
D_{\alpha}D_{p}\Omega$. Then $U_\alpha\in H^{2,2}$. By the Hodge-Riemann
relations we know that,  for any $\alpha$ 
\begin{equation}\label{lasterm1}
(E_{k\alpha i}-U_{\alpha},\bar{E_{k\alpha i}-U_{\alpha}}) \geq 0.
\end{equation}
By~\eqref{fourhodge}, we have
\begin{eqnarray}\label{hodge1}
h_{i\bar j}=\delta_{ij}+\sum_{\alpha}(D_{\alpha}D_{i}\Omega,
\bar{D_{\alpha}D_{j}\Omega}).
\end{eqnarray}
Thus
\begin{align}\label{lasterm2}
\begin{split}
\sum_{\alpha}(U_{\alpha},\bar{U_{\alpha}})
=& \sum_{\alpha}(h^{p\bar q}
\sum_{\beta}(E_{k\beta
i},\bar{D_{\beta}D_{q}\Omega}))(\bar{h^{p_{1}\bar{q_{1}}}
\sum_{\gamma}(E_{k\gamma i},\bar{D_{\gamma}D_{q_{1}}\Omega})}) 
 (D_{\alpha}D_{p}\Omega,\bar{D_{\alpha}D_{p_{1}}\Omega})  \\
=& (h^{p\bar q}
\sum_{\beta}(E_{k\beta i},\bar{D_{\beta}D_{q}\Omega}))(h^{q_{1}\bar{p_{1}}}
\sum_{\gamma}(D_{\gamma}D_{q_{1}}\Omega,\bar{E_{k\gamma i}}))
(h_{p\bar{p_{1}}}-\delta_{pp_{1}})  \\
=& h^{q_1\bar q}\sum_{\beta}(E_{k\beta i},\bar{D_{\beta}D_{q}\Omega})
\sum_{\gamma}(D_{\gamma}D_{q_{1}}\Omega,\bar{E_{k\gamma i}})  
 \\
&-\sum_{p}h^{p\bar q}
\sum_{\beta}(E_{k\beta i},\bar{D_{\beta}D_{q}\Omega})h^{q_{1}\bar{p}}
\sum_{\gamma}(D_{\gamma}D_{q_{1}}\Omega,\bar{E_{k\gamma i}})
 \\
=& h^{p\bar q}\sum_{\beta}(E_{k\beta i},\bar{D_{\beta}D_{q}\Omega})
\sum_{\gamma}(D_{\gamma}D_{p}\Omega,\bar{E_{k\gamma i}}) 
-\sum_{p}\left | h^{p\bar q}\sum_{\beta,q}(E_{k\beta
i},\bar{D_{\beta}D_{q}\Omega}) \right |^{2}
  \\
\leq &  h^{p\bar q}\sum_{\beta}(E_{k\beta i},\bar{D_{\beta}D_{q}\Omega})
\sum_{\gamma}(D_{\gamma}D_{p}\Omega,\bar{E_{k\gamma i}}) \\
=&\sum_{\beta}(E_{k\beta i},\bar{U_{\beta}}),
\end{split}
\end{align}
and
\begin{eqnarray}\label{lasterm3}
\sum_{\alpha}(U_{\alpha},\bar{E_{k\alpha i}})=
h^{p\bar q}
\sum_{\beta}(E_{k\beta i},\bar{D_{\beta}D_{q}\Omega})
\sum_{\alpha}(D_{\alpha}D_{p}\Omega,\bar{E_{k\alpha i}}).
\end{eqnarray}
Combining~\eqref{lasterm1}, \eqref{lasterm2} and
\eqref{lasterm3}  we have
\begin{equation}\label{lasterm4}
 \sum_{\alpha}(E_{k\alpha i}, \bar{E_{k\alpha i}})  
-h^{p\bar q}
\sum_{\alpha}(E_{k\alpha i},\bar{D_{\alpha}D_{q}\Omega})
\sum_{\beta}(D_{\beta}D_{p}\Omega,\bar{E_{k\beta i}})\geq 0.
\end{equation}
Thus the sum of the last two terms in ~\eqref{h2nd3} is nonnegative. 

We shall show that the term
\[
\sum_{\alpha} (D_kD_\alpha D_i\Omega, \bar{D_kD_\alpha D_i\Omega})
\]
is related to the Yukawa coupling of fourfolds.

\begin{definition}
Using the same notations as above, define a holomorphic section of
${\rm Sym}^2\underline{F}^4\otimes (T^*{\mathcal M})^{\otimes 4}$ to be
\begin{equation}\label{Yuka}
\xi_{ijkl}=(\Omega,\pa_i\pa_j\pa_k\pa_l\Omega).
\end{equation}
We call $\xi_{ijkl}$  the Yukawa coupling for  Calabi-Yau fourfolds. 
\end{definition}
Clearly $\xi_{ijkl}$ is symmetric with respect to $i,j,k,l$.

\begin{lemma}\label{y1}
Using the same notations as above, we have
\begin{equation}\label{ly1}
\xi_{ijkl}=-(D_jD_k D_l\Omega, D_i\Omega)=(D_kD_l\Omega, D_jD_i\Omega).
\end{equation}
\end{lemma}

{\bf Proof.} 
The lemma follows from the definition of $T_{k\alpha i}$ and  the first
Hodge-Riemann relation.

\qed

Using the above lemma, we have
\begin{equation}\label{y3}
\sum_{\alpha} (D_kD_\alpha D_i\Omega, \bar{D_kD_\alpha D_i\Omega})
=-\sum_{\alpha, l}|\xi_{ik\alpha l}|^2.
\end{equation}

Combining \eqref{h2nd3}, \eqref{lasterm4} and \eqref{y3} we have 
\begin{align}\label{h2nd4}
\begin{split}
&\widetilde{R}_{i\bar{i}k\bar{k}} \geq 
1+\delta_{ik}-2(D_{k}D_{i}\Omega,\bar{D_{k}D_{i}\Omega})
+\sum_{\alpha,\gamma}(D_{\alpha}D_{i}\Omega,\bar{D_{k}D_{\gamma}\Omega})
(D_{k}D_{\gamma}\Omega, \bar{D_{\alpha}D_{i}\Omega})  \\
&+\sum_{\alpha,\beta}(D_{\alpha}D_{i}\Omega,\bar{D_{\beta}D_{i}\Omega})
(D_{\beta}D_{k}\Omega,\bar{D_{\alpha}D_{k}\Omega}) 
-\sum_{\alpha, l}|\xi_{ik\alpha l}|^2.
\end{split}
\end{align}

The quadratic form $(\,,\,)$ defines an inner product on $H^{2,2}$ by the second
Hodge-Riemann relation. 
Let $\omega_{1},\cdot\cdot\cdot ,\omega_{N}$ be a (real) basis of
$H^{2,2}$ such that $(\omega_{p},{\omega_{q}})=\delta_{pq}$. 
Fix the index $i$. Let $D_{i}D_{\alpha}\Omega=\sum_{p=1}^NA_{\alpha
p}\omega_{p}$ and let 
$D_{k}D_{\beta}\Omega=\sum_{p=1}^N B_{\beta p}\omega_{p}$. By  
Lemma \ref{h22} we have 
\begin{align}\label{h2nd5}
\begin{split}
&\sum_{\alpha,\beta}(D_{\alpha}D_{i}\Omega,\bar{D_{\beta}D_{i}\Omega})
(D_{\beta}D_{k}\Omega,\bar{D_{\alpha}D_{k}\Omega}) 
-\sum_{\alpha\beta}|\xi_{ik\alpha\beta}|^2 \\
=&\sum_{\alpha,\beta}(D_{i}D_{\alpha}\Omega,\bar{D_{i}D_{\beta}\Omega})
(D_{k}D_{\beta}\Omega,\bar{D_{k}D_{\alpha}\Omega})
-\sum_{\alpha,\beta}(D_{k}D_{\beta}\Omega,D_{i}D_{\alpha}\Omega)
(\bar{D_{i}D_{\beta}\Omega},\bar{D_{k}D_{\alpha}\Omega})  \\
=&\sum_{j,l=1}^N\sum_{\alpha,\beta}(A_{\alpha j}\bar{A_{\beta j}}
B_{\beta l}\bar{B_{\alpha l}}-A_{\alpha j}B_{\beta j}
\bar{A_{\beta l}}\bar{B_{\alpha l}}).         
\end{split}
\end{align}
Let $u_{jl}=\sum_{\alpha}A_{\alpha j}\bar{B_{\alpha l}}$. 
From \eqref{h2nd5} we have
\begin{align}\label{h2nd6}
\begin{split}
&\sum_{\alpha,\beta}(D_{\alpha}D_{i}\Omega,\bar{D_{\beta}D_{i}\Omega})
(D_{\beta}D_{k}\Omega,\bar{D_{\alpha}D_{k}\Omega}) 
-\sum_{\alpha,\tau}(D_{k}D_{\tau}\Omega,D_{\alpha}D_{i}\Omega)
(\bar{D_{k}D_{\tau}\Omega},\bar{D_{\alpha}D_{i}\Omega})  \\
=&\sum_{j,l=1}^Nu_{jl}\bar{u_{jl}}-\sum_{j,l=1}^Nu_{jl}\bar{u_{lj}}
=\sum_{j<l}|u_{jl}-u_{lj}|^{2} \geq 0 .
\end{split}
\end{align}
Combining \eqref{h2nd4} and \eqref{h2nd6} we have 
\begin{eqnarray}\label{h2nd7}
\widetilde{R}_{i\bar{i}k\bar{k}} \geq 
1+\delta_{ik}-2(D_{k}D_{i}\Omega,\bar{D_{k}D_{i}\Omega})
+\sum_{\alpha,\gamma}|(D_{\alpha}D_{i}\Omega,\bar{D_{k}D_{\gamma}\Omega})|^{2}
.
\end{eqnarray}
If $i \ne k$, then by~\eqref{h2nd7}
\begin{eqnarray}\label{bisec}
\widetilde{R}_{i\bar{i}k\bar{k}} \geq 
1-2(D_{k}D_{i}\Omega,\bar{D_{k}D_{i}\Omega})
+|(D_{k}D_{i}\Omega,\bar{D_{k}D_{i}\Omega})|^{2} \geq 0.
\end{eqnarray}
 This implies the 
holomorphic bisectional curvature of $\omega_{PH}$ is nonpositive. 

Now we estimate the holomorphic sectional curvature. Let $i=k$. By 
\eqref{h2nd7} we have 
\begin{align}\label{holosec4}
\begin{split}
\widetilde{R}_{i\bar{i}i\bar{i}} \geq &
2-2(D_{i}D_{i}\Omega,\bar{D_{i}D_{i}\Omega})
+\sum_{\alpha,\beta}|(D_{\alpha}D_{i}\Omega,\bar{D_{\beta}D_{i}\Omega})|^{2}
 \\
\geq & 2-2(D_{i}D_{i}\Omega,\bar{D_{i}D_{i}\Omega})
+\sum_{\alpha}|(D_{\alpha}D_{i}\Omega,\bar{D_{\alpha}D_{i}\Omega})|^{2}.
\end{split}
\end{align}
By \eqref{hodge1} we have 
\begin{eqnarray}\label{hodge2}
h_{i\bar i}=1+\sum_{\alpha}
(D_{\alpha}D_{i}\Omega,\bar{D_{\alpha}D_{i}\Omega}).
\end{eqnarray}
Let $a_{\alpha}=(D_{\alpha}D_{i}\Omega, \bar{D_{\alpha}D_{i}\Omega})$ for 
$\alpha \ne i$ and let $a_{i}=(D_{i}D_{i}\Omega,\bar{D_{i}D_{i}\Omega})-1$. 
Clearly they are real numbers by the Hodge-Riemann relations. From 
\eqref{holosec4} and \eqref{hodge2} we have 
\begin{eqnarray*}
\widetilde{R}_{i\bar{i}i\bar{i}} \geq 1+\sum_{\alpha} a_{\alpha}^{2},
\end{eqnarray*}
and 
\begin{eqnarray*}
h_{i\bar i}=2+\sum_{\alpha} a_{\alpha}.
\end{eqnarray*}
Combining the above two inequalities and the following trivial inequality
\begin{eqnarray*}
1+\sum_{\alpha=1}^{m} a_{\alpha}^{2} \geq \frac{1}{m+4}
(2+\sum_{\alpha=1}^{m} a_{\alpha})^{2},
\end{eqnarray*}
we have 
\begin{eqnarray}\label{holosec5}
\widetilde{R}_{i\bar{i}i\bar{i}} \geq \frac{1}{m+4} (h_{i\bar i})^{2}.
\end{eqnarray}
This proved the holomorphic sectional curvature of $\omega_{PH}$ is bounded
above  by a negative constant. Clearly the Ricci curvature is bounded above
by the same negative constant since the
bisectional curvature is nonpositive. 

\qed

\section{Scalar curvature bounds the sectional curvature}\label{S5}
In this section we will prove that the volumes of any subvariety of the moduli
space equipped with the Weil-Petersson metric or the Hodge metric
(Definition~\ref{hodgemetric}) are finite. Also, we will show that the Riemannian
sectional curvature of the Weil-Petersson metric is finite in the $L^{1}$
sense. The key tool we use here is Yau's Schwarz Lemma~\cite{Y3}. The following
version is proved by Royden~\cite{Royden}.
\begin{theorem}
Let 
$M$, $N$ be two \ka manifolds such that $M$ is complete and the
Ricci curvature of $M$ is lowerly bounded and the 
holomorphic sectional curvature of $N$ is
upperly bounded by a negative constant.
Then there is a constant $C$, depending
only on the lower bound of the Ricci
curvature $M$ and the upper bound of the
holomorphic sectional curvature of $N$
such that
\[
\omega_N\leq C\omega_M. 
\]
\end{theorem}

Use the above theorem, we first have 
\begin{theorem}\label{vol}
Let $\mathcal{M}$ be the moduli space of polarized Calabi-Yau $n$-folds. Then 
the volume of any subvariety $M_1$ of $\mathcal{M}$ equipped with the
Weil-Petersson metric or the Hodge metric  is finite. 
\end{theorem}

{\bf Proof.}
Since the moduli space is quasi-projective, after desingularization, we
can assume that $\mathcal{M}=Y \setminus R$ where $Y$ is a compact \ka manifold 
and $R$ is a divisor of normal crossings.
From ~\cite{JY1}, we know that there is a complete metric $\omega_0$ on $\mi$
such that its volume is finite and its Ricci curvature has a lower bound.
Moreover, this metric behaves like the Poincare metric near $R$.
 By Theorem 1.2 in \cite{Lu3} the holomorphic sectional 
curvature of the Hodge metric $\omega_{H}$ is negative away from zero. Let
$i$ be the identity map 
\begin{eqnarray}\label{map}
i: (\mathcal{M},\omega_{0}) \to (\mathcal{M},\omega_{H}),
\end{eqnarray}
which is holomorphic. By Schwarz-Yau lemma~\cite{Y3} 
we have 
\begin{eqnarray}\label{small}
\omega_{H}=i^{\ast}\omega_{H} \leq c\omega_{0}
\end{eqnarray}
for some positive constant $c$. Thus
\begin{eqnarray*}
\int_{{\mathcal M}}\omega_{H}^{m} \leq c_{1}\int_{\mathcal M}\omega_{0}^{m} 
< +\infty.
\end{eqnarray*} 
For any subvariety $M_1$ of the moduli space ${\mathcal M}$, we restrict the 
Hodge metric $\omega_H$ to it. By the Gauss equation, the holomorphic
sectional curvature  on the smooth part of the subvariety $M_1$ is negative away
from zero. Since
$M_1$ is either compact  or quasi-projective, using the same argument for
$\mathcal M$, we proved the volume with respect to the Hodge metric is finite.

By Corollary~\ref{cor42},  up to a constant
\begin{eqnarray}\label{wph}
\omega_{WP} \leq \omega_{H}. 
\end{eqnarray}
So the volume of the Weil-Petersson metric on any subvariety of $\mathcal M$ is
also finite. This  finishes the proof. 

\qed

From the above  theorem we can bound the $L^{1}$ norm of the sectional curvature
of  the Weil-Petersson metric. 
\begin{theorem}\label{curvbdd}
Let $\mathcal{M}$ be the moduli space of polarized Calabi-Yau $n$-folds. 
Then the $L^{1}$ norm of the Riemannian sectional curvature of $\mathcal{M}$
equipped with the Weil-Petersson metric is finite. 
\end{theorem}

{\bf Proof.}
For any point $x_0$ in the moduli space, let $z_{1},\cdot\cdot\cdot,z_{m}$ be 
the local normal coordinates at
$x_0$ with respect to the Weil-Petersson metric. Let $X$ and $Y$ be two real
unit  tangent vectors of $\mathcal{M}$ at $x_0$. Clearly there is a constant $c$
which  is independent of $P$ such that 
\begin{align}\label{securv} 
\begin{split}
|R(X,Y,X,Y)|^{2} \leq c |R_{i\bar j k\bar l}|^{2}
=c\sum_{i,j,k,l}R_{i\bar j k\bar l}R_{j\bar i l\bar k}. 
\end{split}
\end{align}
We make the assumptions~\eqref{ass3},  \eqref{ass1} and \eqref{ass2} at $x_0$
like we did in the proof of Theorem \ref{fourdim}. From the basic fact
\begin{eqnarray*}
|(D_{k}D_{i}\Omega,\bar{D_{l}D_{j}\Omega})|^{2} 
\leq \frac{1}{2}(|(D_{k}D_{i}\Omega,\bar{D_{k}D_{i}\Omega})|^{2}
+|(D_{l}D_{j}\Omega,\bar{D_{l}D_{j}\Omega})|^{2})
\end{eqnarray*}
and  the Strominger formula we have 
\begin{align}\label{securv1}
\begin{split}
&|R(X,Y,X,Y)|^{2} 
\leq c\sum_{i,j,k,l}R_{i\bar j k\bar l}R_{j\bar i l\bar k}\\
=& c\sum_{i,j,k,l}(\delta_{ij}\delta_{kl}+\delta_{il}\delta_{kj}
-(D_{k}D_{i}\Omega,\bar{D_{l}D_{j}\Omega})) 
(\delta_{ij}\delta_{kl}+\delta_{il}\delta_{kj}
-(D_{l}D_{j}\Omega,\bar{D_{k}D_{i}\Omega}))\\ =& c
(2m^{2}+2m-4\sum_{i,k}(D_{k}D_{i}\Omega,\bar{D_{k}D_{i}\Omega})
+\sum_{i,j,k,l}|(D_{k}D_{i}\Omega,\bar{D_{l}D_{j}\Omega})|^{2}) \\
\leq & c_{1} (m+\sum_{i,k}(D_{k}D_{i}\Omega,\bar{D_{k}D_{i}\Omega}))^{2}
=c_{1}(\sum_{i}h_{i\bar i})^{2}
\end{split}
\end{align} 
for some  universal constant $c_{1}$ only depending on $m$. Thus from
\eqref{wph} we have 
\begin{eqnarray*}
\int_{\mathcal{M}}|R(X,Y,X,Y)| \omega_{WP}^{m} \leq \sqrt{c_{1}}
\int_{\mathcal{M}} \sum_{i,j}g^{i\bar j}h_{i\bar j}  \omega_{WP}^{m}
\leq m \sqrt{c_{1}} \int_{\mathcal{M}}\omega_{H}^{m} < +\infty. 
\end{eqnarray*}
This proves that the $L^{1}$ norm of the Riemannian sectional curvature of the 
Weil-Petersson is bounded.

 \qed

For the rest of this section ,we assume that $n=4$.
We will prove that the Riemannian sectional curvature of $\omega_\mu$ 
is bounded by the scalar curvature pointwisely up to a constant. The similar
result  has  been proved in \cite{Lu5} in the case of Calabi-Yau threefolds.
\begin{theorem}\label{uplow}
Let $\mathcal{M}$ be the moduli space of polarized Calabi-Yau fourfolds.  
Then there are positive constants $c_{1}$ 
and $c_{2}$ such that the
Riemannian sectional curvature of the partial 
Hodge metric $\omega_{\mu}$ is bounded by
$c_{1}+c_{2}| \widetilde{R}  |$ 
where $\widetilde{R}$ is the scalar curvature of the partial Hodge metric.
\end{theorem}

{\bf Proof.} 
Fix a point $x_0\in \mathcal{M}$ and let $X$ and $Y$ be two real tangent vectors 
at $x_0$ such that $X$ is perpendicular to $Y$ with respect to $\omega_\mu$. Let
$\xi=X-\sqrt{-1}JX$ and $\eta=Y-\sqrt{-1}JY$ where $J$ is  the complex structure
of $\mathcal{M}$. Clearly $\xi$ is perpendicular to $\eta$  with respect
$\omega_\mu$ too. We make the assumption \eqref{ass3}, \eqref{ass1}, and
\eqref{ass2} like we did in the proof of Theorem 
\ref{fourdim} for the local coordinates $(z_1,\cdots,z_m)$ and the local section
$\Omega$. We can choose a unitary transformation of the coordinates such that
$\xi=a\frac{\partial}{\partial z_{i}}$,
$\eta=b\frac{\partial} {\partial z_{j}}$ with $i \ne j$ for some complex numbers
$a,b$ and the matrix $h_{i\bar j}$ of $\omega_\mu$ is diagonalized with $h_{i\bar
j}=\delta_{ij}\lambda_{i}$. We have 
\begin{eqnarray}\label{sec100}
\widetilde{R}(X,Y,X,Y)=\frac{1}{8}(Re(\widetilde{R}(\xi,\bar\eta,\xi,\bar\eta))
-\widetilde{R}(\xi,\bar\xi,\eta,\bar\eta)). 
\end{eqnarray}
In the following, we will use $||v||^{2}$ to denote the square of the norm of a
complex vector with respect to $\omega_\mu$. 
The second term in the right hand side of the above formula is easy to estimate:
\begin{align}\label{bisec100}
\begin{split}
|\widetilde{R}(\xi,\bar\xi,\eta,\bar\eta)|= & |a|^{2}|b|^{2}
\widetilde{R}_{i\bar i j\bar j} = ||\xi||^{2}||\eta||^{2}
\widetilde{R}_{i\bar i j\bar j}\lambda_{i}^{-1}\lambda_{j}^{-1}\\
\leq & ||\xi||^{2}||\eta||^{2}\sum_{i,j}
\widetilde{R}_{i\bar i j\bar j}\lambda_{i}^{-1}\lambda_{j}^{-1}
=||\xi||^{2}||\eta||^{2}|\widetilde{R}|
\end{split}
\end{align}
since $\widetilde{R}_{i\bar i j\bar j} \geq 0$ for $1\leq i,j \leq m$ by
Theorem~\ref{fourdim}.  By Theorem \ref{curvmu} and
Theorem \ref{fourdim},
\begin{align}\label{fake100}
\begin{split}
\widetilde{R}_{i\bar j i\bar j}=& -2F_{i\bar j i\bar j}
+2\sum_{\alpha,\beta}F_{i\bar j \alpha\bar\beta}F_{i\bar j \beta\bar\alpha}
-\sum_{\alpha\beta}\xi_{ii\alpha\beta}\bar{\xi_{jj\alpha\beta}}\\
& +\big ( \sum_{\alpha}(E_{i\alpha i},\bar{E_{j\alpha j}})
-\sum_{p,\alpha,\beta}\lambda_{p}^{-1}(E_{i\alpha i},\bar{D_{\alpha}D_{p}\Omega})
(D_{\beta}D_{p}\Omega,\bar{E_{j\beta j}})
\big ).
\end{split}
\end{align}
Let $G$ be the vector space spanned by $\{ L_{i} \}$ where 
$i=1,\cdot\cdot\cdot,m$ and 
$L_{i}=(E_{i1i},\cdot\cdot\cdot,E_{imi})$. 
We now define a bilinear form  $((\cdot,\cdot))$ on $G$ by
\begin{eqnarray}\label{braket1}
((L_{i},L_{j}))=\sum_{\alpha}(E_{i\alpha i},\bar{E_{j\alpha j}})
-\sum_{p,\alpha,\beta}\lambda_{p}^{-1}(E_{i\alpha i},\bar{D_{\alpha}D_{p}\Omega})
(D_{\beta}D_{p}\Omega,\bar{E_{j\beta j}}).
\end{eqnarray}
By \eqref{lasterm4},  we 
know that $((L_{i},L_{j}))$ is a Hermitian semi-inner product on $G$. So we have
the following Cauchy inequality:
\begin{eqnarray}\label{fake110}
|((L_{i},L_{j}))| \leq \sqrt{((L_{i},L_{i}))((L_{j},L_{j}))}.
\end{eqnarray} 
However, by the proof of Theorem \ref{fourdim} we know 
\begin{align}\label{fake120}
\begin{split}
((L_{i},L_{i}))= \sum_{\alpha}(E_{i\alpha i},\bar{E_{i\alpha i}})
-\sum_{p,\alpha,\beta}\lambda_{p}^{-1}(E_{i\alpha i},\bar{D_{\alpha}D_{p}\Omega})
(D_{\beta}D_{p}\Omega,\bar{E_{i\beta i}}) 
\leq \ \widetilde{R}_{i\bar i i\bar i} \leq |\widetilde{R}| \lambda_{i}^{2}.
\end{split}
\end{align}
Combining \eqref{braket1}, \eqref{fake110} and \eqref{fake120} we have 
\begin{equation}\label{fake130}
\left | \sum_{\alpha}(E_{i\alpha i},\bar{E_{j\alpha j}})
-\sum_{p,\alpha,\beta}\lambda_{p}^{-1}(E_{i\alpha i},\bar{D_{\alpha}D_{p}\Omega})
(D_{\beta}D_{p}\Omega,\bar{E_{j\beta j}}) \right |
\leq |\widetilde{R}| \lambda_{i}\lambda_{j}.
\end{equation}
Since $(-1)(\cdot,\cdot)$ is a Hermitian inner product on $H^{1,3}$, from 
\eqref{h2nd6} we have 
\begin{equation}\label{fake140}
\sum_{\alpha\beta}\xi_{ii\alpha\beta}\bar{\xi_{jj\alpha\beta}}
\leq  \sum_{\alpha,\beta}|F_{i\bar i \alpha\bar\beta}|
|F_{j\bar j \alpha\bar\beta}|.
\end{equation}
By \eqref{holosec4} we have 
\begin{eqnarray}\label{fake150}
\widetilde{R}_{i\bar i i\bar i} \geq 2-2F_{i\bar i i\bar i}
+ \sum_{\alpha,\beta}|F_{i\bar i \alpha\bar\beta}|^{2} 
\geq \frac{1}{2}\sum_{\alpha,\beta}|F_{i\bar i \alpha\bar\beta}|^{2}.
\end{eqnarray}
So combining \eqref{fake140} and \eqref{fake150} we have 
\begin{align}\label{fake160}
\begin{split}
& \sum_{\alpha\beta}\xi_{ii\alpha\beta}\bar{\xi_{jj\alpha\beta}} \leq 
\sqrt{ \sum_{\alpha,\beta}|F_{i\bar i \alpha\bar\beta}|^2
\sum_{\alpha,\beta}|F_{j\bar j \alpha\bar\beta}|^2}
\leq 2 \sqrt{\widetilde{R}_{i\bar i i\bar i}\widetilde{R}_{j\bar j j\bar j}}
\leq 2 |\widetilde{R}| \lambda_{i}\lambda_{j}.
\end{split}
\end{align}
From \eqref{h2nd7} we have 
\begin{align}\label{fake170}
\begin{split}
\widetilde{R}_{i\bar i j\bar j} \geq 1-2F_{i\bar i j\bar j}
+\sum_{\alpha,\beta}|F_{i\bar\beta \alpha\bar j}|^{2} \geq 
\sum_{\alpha,\beta}|F_{i\bar j \alpha\bar\beta}|^{2}
-|F_{i\bar i j\bar j}|^{2}.
\end{split}
\end{align}
So we have 
\begin{align}\label{fake180}
\begin{split}
|2\sum_{\alpha,\beta}F_{i\bar j \alpha\bar\beta}F_{i\bar j \beta\bar\alpha}| 
\leq \sum_{\alpha,\beta} (|F_{i\bar j \alpha\bar\beta}|^{2}
+|F_{i\bar j \beta\bar\alpha}|^{2})
= 2 \sum_{\alpha,\beta}|F_{i\bar j \alpha\bar\beta}|^{2}
\leq 2\widetilde{R}_{i\bar i j\bar j}+2|F_{i\bar i j\bar j}|^{2}.
\end{split}
\end{align}
From \eqref{fake150}, since $i \ne j$ we have 
\begin{eqnarray}\label{fake181}
|F_{i\bar i j\bar j}| \leq \sqrt{2\widetilde{R}_{i\bar i i\bar i}}
\end{eqnarray}
and 
\begin{eqnarray}\label{fake182}
|F_{i\bar i j\bar j}|=|F_{j\bar j i\bar i}|
 \leq \sqrt{2\widetilde{R}_{j\bar j j\bar j}}.
\end{eqnarray}
Combining \eqref{fake180}, \eqref{fake181} and \eqref{fake182} 
we have 
\begin{align}\label{fake190}
\begin{split}
|2\sum_{\alpha,\beta}F_{i\bar j \alpha\bar\beta}F_{i\bar j \beta\bar\alpha}| 
\leq 2\widetilde{R}_{i\bar i j\bar j}+
4\sqrt{\widetilde{R}_{i\bar i i\bar
i}\widetilde{R}_{j\bar j j\bar j}}
\leq |\widetilde{R}|\lambda_{i}\lambda_{j}+4\sqrt{\widetilde{R}^{2}
\lambda_{i}^{2}\lambda_{j}^{2}}=5|\widetilde{R}|\lambda_{i}\lambda_{j}.
\end{split}
\end{align}
From \eqref{fake150} we also have 
\begin{eqnarray}\label{fake200}
|2F_{i\bar j i\bar j}| \leq 2\sqrt{F_{i\bar i i\bar i}F_{j\bar j j\bar j}}
\leq 2(4\widetilde{R}_{i\bar i i\bar i}\widetilde{R}_{j\bar j j\bar j})
^{\frac{1}{4}} \leq 2\sqrt{2|\widetilde{R}|\lambda_{i}\lambda_{j}}
\leq 1+2|\widetilde{R}|\lambda_{i}\lambda_{j}.
\end{eqnarray}
By the Hodge-Riemann relations we know 
$\lambda_{i}=1+\sum_{\alpha}F_{i\bar i \alpha\bar\alpha} > 1$. 
Combining \eqref{fake100}, \eqref{fake130}, \eqref{fake160}, \eqref{fake190} and
\eqref{fake200} we have 
\begin{align}\label{fake210}
\begin{split}
|Re(\widetilde{R}(\xi,\bar\eta,\xi,\bar\eta))| \leq &|a|^{2}|b|^{2}
|\widetilde{R}_{i\bar j i\bar j}| \leq |a|^{2}|b|^{2}+9|\widetilde{R}|
|a|^{2}|b|^{2}\lambda_{i}\lambda_{j} \\
\leq & ||\xi||^{2}||\eta||^{2}(1+9|\widetilde{R}|).
\end{split}
\end{align}
Combining \eqref{sec100}, \eqref{bisec100} and \eqref{fake210} we have 
\begin{align}\label{final100}
\begin{split}
|\widetilde{R}(X,Y,X,Y)| \leq 
\big (\frac{1}{8}+\frac{9}{8}|\widetilde{R}| \big )
||\xi||^{2}||\eta||^{2}=\big (\frac{1}{2}+\frac{9}{8}|\widetilde{R}| \big
) ||X||^{2}||Y||^{2}.
\end{split}
\end{align}
This finishes the proof. 
\qed

\section{Hodge metrics}\label{4}
Let $X,Y$ be finite dimensional Hermitian vector spaces and let $<,>_X,
<,>_Y$ be the Hermitian inner products of $X$ and $Y$, respectively.
Let
$A,B: X\rightarrow Y$ be linear operators. Then we can define the
natural Hermitian inner product for $A,B$ on the space ${\rm
Hom}(X,Y)$ as follows: let $e_1,\cdots,e_n$ be a unitary basis of
$X$. Then define
\begin{equation}\label{definescalar}
<A,B>=\sum_i<Ae_i, B\bar e_i>. 
\end{equation}

Let $D$ be the classifying space defined in Definition~\ref{def22}.
The complexified tangent bundle $TD\otimes\C$ of $D$ can be realized as
the subbundle of
\begin{equation}\label{tangent}
 TD\otimes\C\subset \oplus_{p+q=n}{\rm Hom} (H^{p,q},
H/H^{p,q}).
\end{equation}
By Lemma~\ref{lem122}, $(\sqrt{-1})^{p-q}Q(\quad ,\quad)$ is the
Hermitian inner product  of $H^{p,q}$, it naturally induced a
Riemannian metric
$h$  on
$TD\otimes \C$ via the above realization.

Define an almost complex structure $J$ on $TD\otimes\C$ as follows:
let $X$ be a local section of $TD\otimes\C$, then
\begin{equation}\label{complex}
JX=
\left\{
\begin{array}{ll}
\sqrt{-1}X,& if\quad X\in\Gamma(\oplus_{p+q=n}{\rm Hom}
(H^{p,q},\oplus_{r< p}H^{r,s}))\\
-\sqrt{-1}X, &{if}\quad X\in\Gamma(\oplus_{p+q=n}{\rm Hom}
(H^{p,q},\oplus_{r>p}H^{r,s}))
\end{array}
\right..
\end{equation}

We have the following
\begin{prop}
The Riemannian metric $h$ is $G$ and $J$ invariant. Furthermore,
$J$ is a $G $ invariant complex structure on $D$. Thus
$h$ defines a Hermitian metric on $D$.
\end{prop}

{\bf Proof.} Let $x,y\in D$ and $\xi x=y$ for some $\xi\in G$. Let
$e_1^p,\cdots,e_{n_p}^p$ be an unitary basis of $H^{p,q}$ at x. Then
$\xi e^p_1,\cdots,\xi e_{n_p}^p$ will be an unitary basis of $H^{p,q}$
at
$y$ by the definition of $G$. If $X$ is a tangent vector at $x$, then
$X$ induced the tangent vector $\tilde X=\xi X\xi^{-1}$ at $y$. Thus
\[
||\tilde X||^2=(\sqrt{-1})^{p-q}\sum_p\sum_i Q(\xi X\xi^{-1}\xi
e_i^p,\xi X\xi^{-1}\xi e_i^p)=||X||^2,
\]
which proves the invariance of $h$ with respect to $G$.

To prove that $h$ is also $J$ invariant, we let $X,Y$ be two
holomorphic vectors at $x$. We just need to prove that $h(X,Y)=0$. Let
$p+q=n$. Suppose
$X$ is nonzero restricting to $H^{p,q}$. As above, let
$e_1^p,\cdots,e^p_{n_p}$ be the unitary basis of $H^{p,q}$. We claim
that
\[
Q(Xe^p_i,Y\bar e^p_i)=0,\quad 1\leq i\leq\dim H^{p,q}.
\]
To see this, assume that a component of  $Xe^p_i\in H^{r,s}$. Then we
have
$r<p$. In order that $Q((Xe_i^p)^{r,s}, Y\bar e_i^p)\neq 0$, we must
have
$s<q$. 
But this is a contradiction because
$p+q=r+s=n$.

It remains to prove that the almost complex structure is integrable.
To prove this, we first observe that the same $J$ defines an almost 
complex structure on $\check D$, the compact dual of $D$. The almost
complex structure on $\check D$ is defined by the pull back of the
complex structure of the flag manifold, which is a  complex
manifold. 

\qed

The main result of this section is:

\begin{theorem}\label{hodge20}
Let ${\mathcal M}$ be a horizontal slice of the classifying space $D$
coming from the moduli space of the polarized Calabi-Yau manifolds.
Then the metric $h$ is the Hodge metric on ${\mathcal M}$.
In particular, the Hodge metric is \ka.
\end{theorem}

The assumption in the theorem can be weakened to the case where
there is a horizontal slice with the Weil-Petersson metric is defined.
See \S ~\ref{WPG} for details. 

{\bf Proof.} Let $D_1$ be the Hermitian symmetric space defined by the
set of subspaces 
\[
H^+=H^{n,0}\oplus H^{n-2,2}\oplus\cdots
\]
in $H$. As in Section~\ref{2}, the natural projection 
\[
p: D\rightarrow D_1
\]
is defined by
\[
\{F^k\}_{k=1,\cdots,n}\mapsto H^+,
\]
which is in general not holomorphic.
Using the same method as above, we defined the unique complex structure
on
$D_1$ by realizing the holomorphic tangent bundle of $D_1$ as the
subbundle of
$
{\rm Hom} (H^+, H/H^+).
$

Let $X\in T^{1,0}{\mathcal M}$ be a holomorphic vector field. 
Then $X$ is horizontal in the sense that $X$ is a section of the bundle
${\rm Hom}(H^{n,0}, H^{n-1,1})\oplus{\rm Hom}(H^{n-1,1},
H^{n-2,2})\oplus\cdots$.
Define
$X_1$ to be an element in
the bundle
such that 
\[
X_1=\left\{
\begin{array}{ll}
X & {\rm restricting ~to } \quad H^+\\
0&  {\rm otherwise}
\end{array}
\right..
\]
and let $X_2=X-X_1$. Then we have
\[
X=X_1+X_2.
\]
Furthermore, we have

\begin{lemma}
According to the complex structure of $D_1$, the vectors fields
$X_1$ and $X_2$ are holomorphic and anti-holomorphic, respectively.
\end{lemma}

{\bf Proof.}
Let 
\[
H=H^+\oplus H^-
\]
Then $X_1$ is a map $H^+\mapsto H^-$, which can be identified as 
 a holomorphic vector fields of $D_1$. $X_2$ can be identified as a
map from $H^-$ to $H^+$. It is the dual map of $H^+\rightarrow H^-$
under the polarization $Q$. Thus can be identified as an
anti-holomorphic vector field.

\qed

{\bf Continuation of the proof of Theorem~\ref{hodge20}.}
From the above argument, we see that under the invariant
\ka metric of
$D_1$
\begin{equation}\label{norm}
||X||^2=||X_1||^2+||X_2||^2.
\end{equation}

If $n$ is an odd number, then  $D_1$ is the
Hermitian symmetric space of third kind. That is
\[
D_1=Sp(n,\R)/U(n).
\]
It can be realized as the subset of $n\times n$ complex matrices
\[
\{Z\in M_n(\C)|I_n-\bar{Z}^TZ>0, Z^T=Z\}.
\]
Its invariant \ka metric can be defined as
\[
\bb \pa\bar\pa\,\log\det (I_n-\bar{Z}^TZ).
\]
If $n$ is an even number, then
\[
D_1=O(m,n,\R)/(O(m)\times O(n)).
\]
There is a natural inclusion  of $D_1$:
\[
D_1\hookrightarrow D_1'=SU(m,n\C)/S(U(m)\times U(n)).
\]
$D_1'$  is the Hermitian symmetric space of first kind, which can
be realized as the subset of
$m\times n$ complex matrices
\[
\{Z\in M_{m,n}(\C)|I_n-\bar{Z}^TZ>0\}.
\]
The invariant \ka metric is defined as
\[
\bb\pa\bar\pa\,\log\det (I_n-\bar{Z}^TZ).
\]
The invariant Riemannian metric pn $D_1$ is the pull back of the 
invariant Hermitian metric on $D_1'$.  

In both cases (of $D_1$ for $n$ odd and of $D_1'$ for $n$ even),
 the
invariant
\ka metrics are defined using the polarization $Q$ as
\[
\bb\pa\bar\pa\log\det Q(\,,\,).
\]

In order to prove the theorem, we just need to prove it at the
original point. At the original point of $D_1$ ($D_1'$, {\it
resp}),
 we
can write the
\ka metric as
\[
\bb dZ_{ij}\otimes d\bar Z_{ij}.
\] 
Let $\{\frac{\pa}{\pa z_\alpha}\}_{\alpha=1,
\cdots,m}$ be  holomorphic
vector fields of ${\mathcal M}$, we have
\[
\sum_{ij}\frac{\pa Z_{ij}}{\pa z_\alpha}\cdot\frac
{\pa\bar Z_{ij}}{\pa z_\alpha}=0,\quad 1\leq\alpha\leq m
\]
The reason for the above equality is that each row of the matrix
$Z_{ij}$ represents an element in some $H^{p,q}$. By the Griffiths
transversality, we have
\[
\frac{\pa Z_{ij}}{\pa z_\alpha}\in H^{p-1,
q+1},
\quad
\frac
{\pa\bar Z_{ij}}{\pa z_\alpha}\in
H^{q-1,p+1}.
\]
The inner product of the above two is the same as $Q(\frac{\pa
Z_{ij}}{\pa z_\alpha},\frac{\pa\bar Z_{ij}}{\pa z_\alpha})$, which is
zero. Similarly, we have
\[
\sum_{ij}\frac{\pa Z_{ij}}{\pa \bar z_\alpha}\cdot\frac
{\pa\bar Z_{ij}}{\pa \bar z_\alpha}=0.
\]
Thus we have
\[
||\frac{\pa}{\pa z_\alpha}||^2
=\sum_{ij} ||\frac{\pa Z_{ij}}{\pa z_\alpha}||^2+\sum_{ij}
||\frac{\pa \bar Z_{ij}}{\pa z_{\alpha}}||^2.
\]
Comparing the above equation with ~\eqref{norm}, we proved that 
$h$ is the Hodge metric. Using the result in~\cite{Lu3}, we know that
the $h$ is \ka.

\qed

Let $\omega_H$ be the \ka form of the Hodge metric $h$. Then we have

\begin{cor}
In the case of $n=3$, then we have
\[
\omega_H=(m+3)\omega_{WP}+{\rm Ric} (\omega_{WP}).
\]
In the case of $n=4$, we have
\[
\omega_H=2(m+2)\omega_{WP}+2{\rm Ric} (\omega_{WP}).
\]
\end{cor}

{\bf Proof.} The case $n=3$ was proved in~\cite{Lu5} using the result
of ~\cite{BG}. In the case of $n=4$, we don't have the similar
result as that of~\cite{BG}. However, using the generalized Strominger
formula~\eqref{stro1}, we have
\[
2(m+2)\omega_{WP}+2{\rm Ric}(\omega_{WP})=2\omega_{WP}
+2g^{k\bar l} F_{i\bar jk\bar l}dz_i\wedge d\bar z_j,
\]
where $F$ is defined in~\eqref{yucawa}. 
Let $X$ be a holomorphic vector on $\mathcal M$. Then by the
identification~\eqref{tangent} and the fact that $X$ is horizontal, 
\[
X\in
{\rm Hom} (H^{4,0}, H^{3,1})\oplus {\rm Hom} (H^{3,1}, H^{2,2})
\oplus {\rm Hom} (H^{2,2}, H^{1,3})\oplus
{\rm Hom} (H^{1,3}, H^{0,4}).
\]

Using (3) of  Lemma~\ref{h31}, we know
that
$\omega_{WP}$ gives the part of the  metric $h$ restricted on the space
$H^{4,0}$. 
Since $D_i\Omega$ gives a basis of the space $H^{3,1}$,
by~\eqref{definescalar}, $h$ restricts to $H^{3,1}$ gives
\[
||X||^2_{H^{3,1}}=\sum_{ij} g^{i\bar j}
(X(D_i\Omega),\bar{X(D_j\Omega)}) (\Omega,\bar\Omega)^{-1}.
\]
In particular, if $X=\frac{\pa}{\pa z_\alpha}$, then 
\[
||X||^2_{H^{3,1}}=g^{i\bar j} F_{i\bar j\alpha\bar\alpha}.
\]
To compute the norm of $X$ on $H^{1,3}$ and $H^{0,4}$, we use the
duality as follows:
let $w_1,\cdots,w_N$ be a (real) orthonormal basis of $H^{2,2}$. Then
we have
\[
||X||^2_{H^{2,2}}=\sum_{\alpha=1}^N (X(e_\alpha), \bar{X(e_\alpha)})
=\sum_{ij}\sum_{\alpha=1}^N g^{i\bar j}(\Omega,\bar\Omega)^{-1}
(X(e_\alpha), {D_i\Omega}) (\bar{D_j\Omega}, \bar{X(e_\alpha)}).
\]
Since the polarization is invariant infinitesmally under $X$, the
above is equal to 
\[
||X||^2_{H^{2,2}}=
\sum_{ij}\sum_{\alpha=1}^N g^{i\bar j}(\Omega,\bar\Omega)^{-1}
(e_\alpha, D_\alpha{D_i\Omega}) (\bar{D_\alpha D_j\Omega},
{e_\alpha})=\sum_{ij}g^{i\bar j}(\Omega,\bar\Omega)^{-1}
(D_\alpha D_i\Omega, \bar{D_\alpha D_j\Omega}).
\]
Thus the norm restricted to $H^{2,2}$ is the same as that on
$H^{3,1}$. Using the same method, we can prove that the norm of $X$
on $H^{1,3}$ is given by the Weil-Petersson metric. The corollary
follows from 
\[
||X||_h^2=||X||^2_{H^{4,0}}+||X||^2_{H^{3,1}}+||X||^2_{H^{2,2}}+
||X||^2_{H^{1,3}}.
\]

\qed

\begin{cor}[{cf.~\cite{Lu3}}]\label{cor42}
Up to a constant, 
the Weil-Petersson metric and its Ricci curvature are less
than or equal to the Hodge metric.
\end{cor}

{\bf Proof.} This is an easy consequence
of Theorem~\ref{Str}.

\qed

\begin{rem}
It was a well known theorem of ~\cite{GS} that the holomorphic
sectional curvature with respect to the Hermitian connection at the
horizontal direction is negative away from zero. In the previous
section, we give an explicit formula proving that the holomorphic
bisectional curvature on the horizontal slice is nonpositive and the
holomorphic sectional curvature is negative away from zero
in the case of Calabi-Yau fourfolds.
Since the
Hodge metric is \ka, the connection is also the Levi-Civita connection.
\end{rem}

\begin{rem}
One of 
the most confusing part of the theory of the Hodge metric is that the
projection in Definition~\ref{def24} is, in
general, not holomorphic. 
This is of course true if $D_1$ is not a Hermitian symmetric
space. Even if $D_1$ is a Hermitian symmetric space,
the projection is Definition~\ref{def24} is in general
not holomorphic. However, in this case,
there is a unique complex structure on $D$
that will make the projection holomorphic
and thus make the manifold $D$ homogeneous
K\"ahler. $D$ is in general not homogeneous
K\"ahler, thus the invariant Hermitian
metric can't be a \ka metric.

Take a closer look of the above phenomena~\footnote{This was pointed to the authors by Professor A.
Todorov.}. Let $D=G/V$ as in
Section~\ref{2}. Consider the isotropy representation of the compact
group
$V$ in
$T_0(D)$, the tangent space of $D$ at the original
point. If $V$ is the maximal compact subgroup of $G$, then the
representation is irreducible and thus there is only one invariant
almost complex structure. In general the group representation is not
irreducible. Thus there are $2^N$ different almost complex structures
on $D$ where $N$ is the number of irreducible components of the
representation.
\end{rem}


\section{The curvature of the Hodge metric 
in dimension 1}\label{jjj} 
In this section, we prove that in the one dimensional case,
the  curvature of the Hodge metric is bounded near the boundary points
with infinite Hodge distance. We will consider the $n$-dimensional
case in the next paper~\cite{LS-2}.

Our starting point  is the relation between the completeness
of the metrics and the limiting Hodge structures. Such a relation
was first drawn by C. Wang. In his paper~\cite{Wong}, among the 
other results, Wang proved the following

\begin{theorem}
Let $\Delta^*$ be the one dimensional parameter space of
a family of polarized Calabi-Yau manifolds. Then the necessary
and sufficient condition for the Weil-Petersson mertric
to be complete is $NA_0\neq 0$, where $N$ is the nilpotent operator
in~\eqref{asa}
of $\Delta^*$ and $A_0$ is defined in~\eqref{Aseries}.
\end{theorem}

As above, let $\Delta^*$ be the one dimensional parameter space of
a family of polarized Calabi-Yau manifolds. Let $\Omega$
be the section of the first Hodge bundle $F^n$. Then by
the Nilpotent Orbit theorem of Schmid (Theorem~\ref{tnot}),
after
a possible base change,  we have
\begin{equation}\label{asa}
\Omega=e^{\bb N\log\frac 1z}A(z),
\end{equation}
where $N$ is the nilpotent operator, $N^{n+1}=0$ for 
$n$ the dimension of the Calabi-Yau manifolds, 
and 
\begin{equation}\label{Aseries}
A(z)=A_0+A_1z+\cdots
\end{equation}
is a
vector valued convergent power series with the
convergent radius $\delta>0$.
(see Section~\ref{2} for details). Let  
\[
f_{k,l}(z)=z^k(\log\frac 1z)^l,
\]
for any $k,l\geq 0$. Then we can write $\Omega$ as the
convergent series
\begin{equation}\label{bseries}
\Omega=\sum_{k,l}A_{k,l} z^k(\log\frac 1z)^l=\sum_{k,l}
A_{k,l}f_{k,l}.
\end{equation}

Define $\deg f_{k,l}=k-\frac{l}{n+1}$. Then
we have the following lemma:

\begin{lemma}\label{sp1}
The convergence of ~\eqref{bseries} is in the
$C^\infty$ sense. Furthermore, we have
\begin{equation}\label{esti}
||\Omega-\sum_{\deg
f_{k,l}\leq\mu}A_{k,l}f_{k,l}||_{C^s}\leq C
r^{k_0-s}(\log\frac 1r)^{l_0},
\end{equation}
where $r=|z|$, $k_0,l_0$ are the unique pair of nonnegative
integers such that
$l_0\leq n$,
$k_0-\frac{l_0}{n+1}>\mu$  and for any pair of
integers $k',l'$ with 
$k'-\frac{l'}{n+1}>\mu$ we have
$k'-\frac{l'}{n+1}\geq k_0-\frac{l_0}{n+1}$. 
$C$ is  a constant depending only on $k_0, l_0, \mu$ and
$\Omega$.
\end{lemma}

{\bf Proof.}  From ~\eqref{Aseries}, we have $|A_k|\leq
(\frac\delta 2)^{-k}$ and thus $|A_{k,l}|\leq
(\frac\delta 4)^{-k-1}$ for small $\delta$ and large
$k$. Thus we know that
\[
\sum_{k,l}|A_{k,l}f_{k,l}|\leq\sum(\frac\delta 4)^{-k-1}
r^{k}(\log\frac 1r)^l<+\infty,
\]
and thus the convergence  in ~\eqref{bseries} is
uniform for $r<\delta/4$. To prove that the convergence is $C^s$ for any
$s\geq 1$, we observe that
\[
f_{k,l}'=kf_{k-1,l}-lf_{k-1,l-1}.
\]
Thus we have
\begin{equation}\label{esti2}
\sum_{k,l}
|A_{k,l}kf_{k-1,l}|+|A_{k,l}lf_{k-1,l-1}|\leq\sum_{k,l}
(\frac\delta 4)^{-k-1}r^{k-1}(\log\frac
1r)^l(k+n)<+\infty,
\end{equation}
for $r<\delta/4$. Thus the convergence is $C^1$. Using
mathematical induction, the convergence is in fact in
the $C^k$ sense. To get the quantitative
result~\eqref{esti}, we just observe that
$(f_{k,l})^{(s)}$ is a linear combination of
$f_{k-s,l},\cdots, f_{k-s, l-s}$ with the coefficients
not more than $(2(|k|+|l|))^s$. An
inequality like ~\eqref{esti2} gives the requires
estimates.

\qed

Having finished the convergence of the series, 
we prove the following

\begin{theorem}\label{bound}
Assume the moduli space $\mathcal M$ of polarized 
Calabi-Yau threefolds is one dimensional. 
If $\Delta^*$ is a holomorphic chart of $\mathcal M$ such 
that $\Delta^*$ 
is complete at $0$ with respect to the Hodge
metric, then the Gauss curvature of the Hodge metric is
bounded.\footnote{The referee pointed out that the result is also true
for partial Hodge metric.} 
\end{theorem}

To prove the theorem,
we first assume that $NA_0\neq 0$. Under this
assumption, we have the following

\begin{lemma}\label{ll11}
If $NA_0\neq 0$, then the expression 
\[
(e^{\bb N\log\frac 1z}A_0,\bar{e^{\bb N\log\frac
1z}A_0}),
\]
is a non-constant polynomial of $\log \frac 1r$, where
$r=|z|$.
\end{lemma}

{\bf Proof.}
By the definition of the operator $N$, we know that
$N$ is an element of the Lie algebra of the Lie group
$G_R$. Thus we know that the above expression is equal
to 
\[
(e^{\bb N\log\frac{1}{r^2}}A_0,\bar{A_0}).
\]

If the above expression is a constant, we then would
have
\begin{equation}\label{use1}
(N^lA_0,\bar{A_0})=0
\end{equation}
for any positive integer $l$. Thus we would have
\[
(\pa_z e^{\bb N\log\frac 1z}A_0, \bar{e^{\bb N\log\frac
1z}A_0})=0.
\]
Since by the assumption, $\pa_z e^{\bb N\log\frac
1z}A_0\neq 0$, the nilpotent orbit theorem implies that
\[
(\pa_z e^{\bb N\log\frac 1z}A_0,\bar{\pa_z e^{\bb
N\log\frac 1z}A_0})<0,
\]
which  is a contradiction.

\qed

In what follows we use $l$ to denote the degree of the polynomial in 
the above lemma.

\begin{cor}\label{cor111} 
If $NA_0\neq 0$, then 
\[
r^{s+2}(\log
\frac
1r)^3||\omega_{WP}-\frac l4\cdot\frac{1}{r^2(\log\frac
1r)^2}dz\wedge d\bar z||_{C^s}\leq C,
\]
for any  integer $s\geq 0$, where $C$ is a constant
depending only on $s$, $n$ and the convergence radius
$\delta$.
\end{cor}

{\bf Proof.}
For any mononomials of
the form
$z^t\bar z^s(\log\frac 1r)^l$, with integers $t, s, l$,
we define the degree of it to be $t+s-l/(n+1)$.
We write
\[
(\Omega,\bar\Omega)
=c(\log\frac 1r)^l+R_0(\log\frac 1r)+
\tilde R(z,\bar z,\log\frac 1r),
\]
where
$c(\log\frac 1z)^l$ is the highest order term of the 
polynomial 
in Lemma~\ref{ll11}, 
$R_0$ is the polynomial of $\log\frac 1r$ of  degree less
than or equal to $l-1$ and $\tilde R(z,\bar z,\log\frac
1r)$ contains the terms with degree at least positive.
From Lemma~\ref{sp1}, the above series converges in the
sense of $C^\infty$. The corollary follows from the
fact that
\[
\omega_{WP}-\frac l4\cdot\frac{1}{r^2(\log\frac
1r)^2}dz\wedge d\bar z
=-\pa\bar\pa\log\frac{(\Omega,\bar\Omega)}{(\log\frac
1r)^l}.
\]

\qed

Now we assume that $NA_0=0$. We normalize $(\Omega,\bar\Omega)$ 
such that $(A_{0},\bar{A_{0}})=1$. Then we
have the  following expansion 
\begin{equation}\label{use2}
\log\, (\Omega,\bar\Omega)=P+\bar P+f(z,\bar
z)(\log\frac 1r)^l+R(z,\bar z),
\end{equation}
where $f(z,\bar z)\not\equiv 0$ is a homogeneous
polynomial of degree $2k$
\footnote{We shall prove that the degree of the 
polynomial is actually an even number.}
, and $P$ is a polynomial of
$z$ of degree less than or equal to  $2k-1$ but
no less than 1  and
$R(z,\bar z)$ is a series of mononomials of degree
great than
$2k-l/(n+1)$. 
In the expansion, we allow that $l=0$. But if $l=0$, we
assume that $f(z,\bar z)$ is not of the form of
$c(z^{2k}+\bar z^{2k})$, otherwise, we can include
$f(z,\bar z)$ in $P+\bar P$. By Lemma~\ref{sp1}, the
expansion is convergent in the
$C^\infty$ sense. We have the following observation

\begin{lemma}\label{nonzero}
If $l\geq 1$, then
there are no $z^{2k}$ or $\bar z^{2k}$ terms  in the
polynomial $f(z,\bar z)$. In particular, 
\[
\pa_z\bar\pa_z f(z,\bar z)\not\equiv 0.
\]
\end{lemma}

{\bf Proof.} From~\eqref{use2}, we have the following
expansion 
\[
(\Omega,\bar\Omega)=1
+P+\bar P+f(z,\bar z)(\log \frac 1r)^l+\cdots,
\]
where the terms in $\cdots$ are the terms of degree 
at least
$2k-(l-1)/(n+1)$, or the terms without $\log\frac 1r$.
If there is a nonzero $z^{2k}$ term in  $f(z,\bar z)$,
we must have
\[
(N^lA_{2k}, \bar A_{0})\neq 0,
\]
which is not possible because of the assumption
$NA_0= 0$. Thus there are no
$z^{2k}$ or $\bar z^{2k}$ terms. Since $f$ is not
identically zero, this implies
\[
\pa_z\bar\pa_z f(z,\bar z)\not\equiv 0.
\]

\qed

By~\eqref{use2} and the $C^\infty$ convergence, we have
\begin{equation}\label{poi12}
\lambda=-\pa_z\bar\pa_z f(z,\bar z) (\log\frac
1r)^l+R(z,\bar z,\log\frac 1r),
\end{equation}
where $\lambda dz\otimes d\bar z$ defines the Weil-Petersson metric
and where $R(z,\bar z,\log\frac 1r)$ contains terms of
degree no less than $2k-2-\frac{l-1}{n+1}$. Since
$\lambda>0$, we must have
\[
-\pa_z\bar\pa_z f(z,\bar z) \geq 0.
\]
Thus $2k$ is an even number, otherwise the integral of
the above expression along the unit circle would be
zero, contradicting to Lemma~\ref{nonzero}. So $k$ is
actually an integer.

\begin{lemma}\label{use3}
Using the same notations as above,
we have
\[
f(z,\bar z)=c r^{2k},
\]
for some constant $c$.
\end{lemma}

{\bf Proof.} By Corollary~\ref{cor42}, we have, up to a constant
\[
h\geq-\pa_z\bar\pa_z\log\lambda,
\]
where $h dz\otimes d\bar z$ defines the Hodge metric.
By the Schwartz-Yau Lemma, we have
\begin{equation}\label{YauS}
-\pa_z\bar\pa_z\log\lambda\leq h\leq
 \frac{1}{r^2(\log\frac 1r)^2}
\end{equation}
up to a constant.
However, at a point where $\pa_z\bar\pa_z f\neq 0$, we
have
\[
\log\lambda=\log (-\pa_z\bar\pa_z f(z,\bar z))+
l\log(\log\frac 1r)+\log(1+\frac{R(z,\bar z,\log\frac
1r)}{-\pa_z\bar\pa_z f(z,\bar z) (\log\frac
1r)^l}).
\]
Using the same method as in the proof of
Corollary~\ref{cor111}, we have
\[
\pa_z\bar\pa_z \log(1+\frac{R(z,\bar z,\log\frac
1r)}{-\pa_z\bar\pa_z f(z,\bar z) (\log\frac
1r)^l})=O(\frac{1}{r^2(\log\frac 1r)^3}).
\]
 
Using~\eqref{YauS}, we have
\[
\pa_z\bar\pa_z\log (-\pa_z\bar\pa_z f(z,\bar z))\equiv0,
\]
otherwise it could have been of the order $r^{-2}$,
which is a contradiction to~\eqref{YauS}.
An elementary argument using Lemma~\ref{nonzero} shows
the
$f(z,\bar z)$ must be of the form stated in the lemma.

\qed

{\bf Proof of Theorem~\ref{bound}.} First we compute the scalar 
curvature of the Hodg metric. We use the
same notation as in the previous sections. Let $\lambda$ be the 
Weil-Petersson metric and let $h$ be the Hodge metric. 
Let  
\begin{align*}
& K=-\log (\Omega,\bar\Omega);\\
&\Gamma_{11}^1=\frac{\pa\log\lambda}{\pa z}=\lambda^{-1}\pa_z\lambda;\\
&K_1=-\pa_z\log (\Omega,\bar\Omega);\\
& F_{111} =(\Omega,\pa_z\pa_z\pa_z\Omega);\\
&F_{1111}=\pa_1 F_{111}-3\Gamma_{11}^1
F_{111}  +2K_1 F_{111};\\
&A=\lambda^{-2}e^{2K}|F_{111}|^2.
\end{align*}
Let $R_{1\bar 11\bar 1}$ be the curvature of the Weil-Petersson metric 
and let $\tilde R_{1\bar 11\bar 1}$ be the curvature of the Hodge metric. 

Since $\mathcal M$ is the moduli space of polarized 
Calabi-Yau threefolds and $\mathcal M$ is one dimensional, from the 
Strominger formula we have 
\[
R_{1\bar 11\bar 1}=2\lambda^2-\lambda^{-1}e^{2K}|F_{111}|^2.
\]
So the Ricci curvature of the Weil-Petersson metric is 
\[
Ric(\lambda)=-\pa_z\pa_{\bar z}\log\lambda=-\lambda^{-1}
R_{1\bar 11\bar 1}=-2\lambda+\lambda^{-2}e^{2K}|F_{111}|^2
=-2\lambda+A.
\]
This implies 
\[
h=(m+3)\lambda+Ric(\lambda)=4\lambda+(-2\lambda+A)
=2\lambda+A=\lambda(2+\lambda^{-3}e^{2K}|F_{111}|^2).
\]
So we have 
\begin{align}\label{add500}
\begin{split}
\pa_z h=&\pa_z\lambda(2+\lambda^{-3}e^{2K}|F_{111}|^2)\\
&+\lambda
[-3\lambda^{-4}\pa_z\lambda e^{2K}|F_{111}|^2
-2\lambda^{-3}(\Omega,\bar\Omega)^{-3}(\pa_z\Omega,\bar\Omega)
|F_{111}|^2
+\lambda^{-3}e^{2K}\pa_z F_{111}\bar{F_{111}}]\\
=&h\lambda^{-1}\pa_z\lambda+\lambda^{-2}e^{2K}\bar{F_{111}}
(-3\Gamma_{11}^1 F_{111}+2K_1 F_{111}+\pa_z F_{111})\\
=&h\Gamma_{11}^1+\lambda^{-2}e^{2K}\bar{F_{111}}F_{1111}.
\end{split}
\end{align}
Similarly, we have 
\[
\pa_{\bar z}h=
h\bar{\Gamma_{11}^1}+\lambda^{-2}e^{2K}F_{111}\bar{F_{1111}}.
\]
So the curvature of the Hodge metric is 
\begin{align}\label{add520}
\begin{split}
\tilde R_{1\bar 11\bar 1}=&\pa_z\pa_{\bar z}h-h^{-1}\pa_z h\pa_{\bar z}h
=\pa_{\bar z}
(h\Gamma_{11}^1+\lambda^{-2}e^{2K}\bar{F_{111}}F_{1111})
-h^{-1}\pa_z h\pa_{\bar z}h\\
=&\pa_{\bar z}h\Gamma_{11}^1+h\pa_{\bar z}\Gamma_{11}^1
-2\lambda^{-3}\pa_{\bar z}\lambda e^{2K}\bar{F_{111}}F_{1111}
-2\lambda^{-2}(\Omega,\bar\Omega)^{-3}(\Omega,\bar{\pa_z\Omega})
\bar{F_{111}}F_{1111}\\
&+\lambda^{-2} e^{2K}\bar{\pa_z F_{111}}F_{1111}
+\lambda^{-2} e^{2K}\bar{F_{111}}\pa_{\bar z}F_{1111}
-h^{-1}\pa_z h\pa_{\bar z}h\\
=&(h\bar{\Gamma_{11}^1}+\lambda^{-2}e^{2K}F_{111}\bar{F_{1111}})
\Gamma_{11}^1+h\pa_{\bar z}\Gamma_{11}^1
-3\lambda^{-2}\Gamma_{11}^1e^{2K}\bar{F_{111}}F_{1111}\\
&+\lambda^{-2}\Gamma_{11}^1e^{2K}\bar{F_{111}}F_{1111}
+2\lambda^{-2}e^{2K}K_1\bar{F_{111}}F_{1111}
+\lambda^{-2}e^{2K}\bar{\pa_z F_{111}}F_{1111}\\
&+\lambda^{-2} e^{2K}\bar{F_{111}}
\pa_{\bar z}(\pa_z F_{111}-3\Gamma_{11}^1 F_{111}+2K_1 F_{111})
-h^{-1}\pa_z h\pa_{\bar z}h.
\end{split}
\end{align}
Using $\pa_{\bar z}\Gamma_{11}^1=\pa_z\pa_{\bar z}\log\lambda
=-Ric(\lambda)=2\lambda-A$, from the above formula we have 
\begin{align}\label{add540}
\begin{split}
\tilde R_{1\bar 11\bar 1}=&h|\Gamma_{11}^1|^2
+\Gamma_{11}^1\lambda^{-2}e^{2K}F_{111}\bar{F_{1111}}
+(2\lambda+A)(2\lambda-A)+\lambda^{-2}e^{2K}|F_{1111}|^2\\
&+\bar{\Gamma_{11}^1}\lambda^{-2}e^{2K}\bar{F_{111}}F_{1111}
+\lambda^{-2}e^{2K}\bar{F_{111}}(-3(2\lambda-A)+2\lambda)F_{111}\\
&-h^{-1}(h\Gamma_{11}^1+\lambda^{-2}e^{2K}\bar{F_{111}}F_{1111})
(h\bar{\Gamma_{11}^1}+\lambda^{-2}e^{2K}F_{111}\bar{F_{1111}})\\
=&4\lambda^2-A^2+\lambda^{-2}e^{2K}|F_{1111}|^2+A(3A-4\lambda)
-h^{-1}\lambda^{-4}e^{4K}|F_{111}|^2|F_{1111}|^2\\
=&4\lambda^2-4\lambda A+2A^2+\lambda^{-2}e^{2K}|F_{1111}|^2(1-h^{-1}A)\\
=&4\lambda^2-4\lambda A+2A^2+\lambda^{-2}e^{2K}|F_{1111}|^2
(2\lambda h^{-1})\\
=&4\lambda^2-4\lambda A+2A^2+2\lambda^{-1}e^{2K}|F_{1111}|^2 h^{-1}.
\end{split}
\end{align}

The scalar
curvature of the Hodge metric is given by
\begin{align}\label{kk1}
\begin{split}
\rho=&-h^{-2}\tilde R_{1\bar 11\bar 1}=
-\frac{4\lambda^2-4\lambda A+2A^2}{(2\lambda+A)^2}
-\frac{2\lambda^{-1}e^{2K}|F_{1111}|^2}{(2\lambda+A)^2}\\
=&
-\frac
{4-4e^{2K}\lambda^{-3}|F_{111}|^2+2e^{4K}\lambda^{-6}
|F_{111}|^4}
{(2+e^{2K}\lambda^{-3}|F_{111}|^2)^2}
-
\frac{2e^{2K}\lambda^{-4}|F_{1111}|^2}
{(2+e^{2K}\lambda^{-3}|F_{111}|^2)^3}.
\end{split}
\end{align}

Apparently, the first term on the right hand side of
~\eqref{kk1} is bounded. Thus in order to prove the
theorem, we just need to bound the second term of the
right hand side of ~\eqref{kk1}.

Case 1. $NA_0\neq 0$. In this case, by
Corollary~\ref{cor111}, we have
\begin{equation}\label{lam}
\lambda\sim\frac{1}{r^2(\log\frac 1r)^2}.
\end{equation}

For the Yukawa coupling $F_{111}$, we always have
$F_{111}=O(\frac{1}{r^3})$. 
If $|F_{111}|=O(\frac{1}{r^2})$, then $|F_{1111}|=O(
\frac{1}{r^3})$. Thus
\[
2e^{2K}\lambda^{-4}|F_{1111}|^2\rightarrow 0,
\]
and is bounded. If $F_{111}\sim\frac{1}{z^3}$, then we
have the following asymptotic computations
\begin{align*}
& \pa_1 F_{111}\sim\frac{-3}{z^4};\\
& \Gamma_{11}^1F_{111}\sim\frac{-1}{z^4};\\
& |K_1F_{111}|\leq C\frac{1}{r^4\log \frac 1r}.
\end{align*}
Thus we have
\begin{equation}\label{F111}
|F_{1111}|\leq C\frac{1}{r^4\log\frac 1r}.
\end{equation}
Using the facts that $F_{111}\sim\frac{1}{z^3}$, we have
\begin{equation}\label{e2k}
e^{2K}\lambda^{-3}\sim r^6.
\end{equation}
Using ~\eqref{lam},~\eqref{F111} and~\eqref{e2k}, we
proved that in this case the curvature is bounded.

Case 2. $NA_0=0$. 
In this case, by ~\eqref{poi12} and Lemma ~\ref{use3},
\[
\lambda=-ck^2 r^{2(k-1)} (\log\frac
1r)^l+R(z,\bar z,\log\frac 1r),
\]
where $R(z,\bar z,\log\frac 1r)$ contains terms of order
at least $2(k-1)-\frac{l-1}{n+1}$. 
We claim that $l\geq 1$, otherwise, by the above
equation, we would have that the Hodge metric, as the
linear combination of the Weil-Petersson metric and
its Ricci curvature, satisfying
\[
h\leq\frac{(\log\frac 1r)^s}{r},
\]
for some positive integer $s$,
and thus is incomplete. A
straightforward computation gives
\[
e^{2K}\lambda^{-2}|F_{111}|^2\sim\frac{1}{r^2(\log
\frac
1r)^2}.
\]
This implies that
\[
F_{111}\sim{z^{2k-3}},
\]
and by using the same argument as we did in Case 1, 
we have
\[
|F_{1111}|\leq{C}{r^{2k-4}}.
\]
Thus we have
\[
\frac{2e^{2K}\lambda^{-4}|F_{1111}|^2}
{(2+e^{2K}\lambda^{-3}|F_{111}|^2)^3}\leq
2e^{-4K}\lambda^5|F_{1111}|^2/|F_{111}|^6,
\]
and it is bounded.

\qed

\section{The Weil-Petersson Geometry}\label{WPG}

 By a classical result of Wolpert~\cite{Wolpert}, the
curvature of the Weil-Petersson metric on 
Teichm\"uller space is nonpositive.  However, the curvature
of the Weil-Petersson metric on the moduli space of Calabi-Yau manifolds
doesn't have such  a good property~\footnote{In fact, physicists
found that the curvature of the Weil-Petersson metric on
certain moduli space can either be positive or
negative~\cite[page 65]{COGP}.}. The bad curvature property
makes it difficult to do geometric analysis on the moduli
space. In order to overcome this difficulty, in ~\cite{Lu3}
and ~\cite{Lu5}, the first author introduced a new \ka metric
called Hodge metric. On one side, the  holomorphic
bisectional curvature of the Hodge metric is nonpositive, on
the other side, up to a constant, the Weil-Petersson metric
is smaller than the Hodge metric. Thus one  can use the
Hodge metric to do the similar geometric analysis as that on 
Teichm\"uller  space and then translate the results back in the language
of  the Weil-Petersson metric.

In the proof of the nonpositivity  of the
curvature of  Hodge metrics (cf. ~\cite{Lu3},~\cite{Lu5}
and Theorem~\ref{fourdim}), we don't need the assumption
that the manifold is the moduli space of Calabi-Yau
manifolds. All we need  is the fact that the manifold is a
horizontal slice and there is a Weil-Petersson metric on
it. In fact, the existence of the Weil-Petersson metric gives
severe
restrictions on the variation of the Hodge structures. 
These kinds of restrictions haven't been studied
comprehensively.

Lemma~\ref{nonzero} and ~\ref{use3} are good examples of
how the existence of the Weil-Petersson metric affects the
variation of the Hodge structures at infinity of the horizontal slices. In
fact, using the notations in \S ~\ref{jjj}, Lemma~\ref{use3} implies the
following

\begin{prop}
Let $k,l$ be defined in ~\eqref{use2}. Then if $l\geq 1$, we
have
\[
(N^lA_p,\bar A_q)=0
\]
for any $p+q=2k$ but $p\neq q$, where the vectors $A_p$ are
defined in ~\eqref{Aseries}.
\end{prop}

Besides the case $p=0$, it is rather difficult to prove the
above result without using the Schwarz-Yau inequality. We
believe that there are more properties of this kind. Because of this, we
defined the following 
concept of the Weil-Petersson geometry and would like to 
study the  properties in a systematic way.

\begin{definition} The Weil-Petersson geometry contains a \ka
orbifold
$M$ with the orbifold metric
$\omega_{WP}$ such that:
\begin{enumerate}
\item Let $\tilde M$ be the universal covering space of $M$. Then
there is a natural immersion $\tilde M\rightarrow D$ from $M$ to the
classifying space $D$ (cf. ~\cite{Gr}) such that $M$ is a
horizontal slice of
$D$. In this way, we can also endowed the Hodge bundles
$F^1,\cdots,F^n$ to $M$ where $F^n$ is a line bundle;
\item $\omega_{WP}$ is the curvature of the bundle $F^n$. It
is positive definite and thus defines a \ka metric in
${M}$ and is called the Weil-Petersson metric;
\item $M$ is quasi-projective and $F^n$ is an ample line
bundle of $M$. 
The compactification is called Viehweg
compactification~\cite[page 21, Theorem 1.13]{V2}. The Hodge
bundles
$F^1,\cdots, F^n$ extend to the compactification $\bar M$ of $M$~
\footnote{This follows from Schmid's Nilpotent 
Orbit Theorem~\cite{Schmid}.};
\item After passing to a finite covering
and after desingularization, in a neighborhood of the infinity,
$M$ can be written as
\[
\Delta^{n-k}\times (\Delta^*)^{k},
\]
where $\Delta$ is the unit disk and $\Delta^*$ is the
punctured unit disk. Let $\Omega$ be a local section of
$F^n$ in the neighborhood, then locally, $\Omega$ can be
(multi-valuedly) written as
\[
\Omega=e^{\sqrt{-1}(N_1\log\frac{1}{z_1}+\cdots+N_k\log\frac{1}{z_k}
)}A(z_1,\cdots, z_n),  
\]
where $N_1,\cdots,N_k$ are nilpotent operators
and $A$ is a vector valued holomorphic function of
$z_1,\cdots,z_n$.
\end{enumerate}
\end{definition}

\begin{rem}
The first property of above is basically the Griffiths transversality
~\cite{Gr}. The second property is a theorem of
Tian~\cite{T1}. The third one is the compactification theorem of
Viehweg~\cite{V2} and the fourth property is the Nilpotent Orbit
theorem of Schmid~\cite{Schmid}.
\end{rem}

The theorems in this paper are true for abstract 
Weil-Petersson geometry defined above. A further study if
the Weil-Petersson geometry will be the project of 
future study. In particular, we wish to define a natural
metric which is a modification of the Hodge metric at
infinity similar to that of McMullen's~\cite{McM} in the
case of Teichm\"uller space. It would be interesting if we can do so
in the category of the Weil-Petersson Geometry.

\section{Appendix}
In this appendix we prove Theorem \ref{curvmu}. 
As before, the subscripts $i,j,\cdots$ are all ranging from $1$ to $m$, unless
otherwise noted.

{\bf Proof of Theorem~\ref{curvmu}.} 
By definition $\omega_{\mu}=\mu\omega_{WP}+Ric(\omega_{WP})$, since 
the Weil-Petersson metric is K\"ahler, we know $\omega_{\mu}$ is $d$-closed. 
From the Strominger formula \eqref{stro1}, we know the Ricci tensor of the 
Weil-Petersson metric is 
\begin{equation}\label{Ric}
R_{i\bar{j}} = -(m+1)g_{i\bar{j}}+g^{k\bar{l}}F_{i\bar jk\bar l}.
\end{equation}
Thus we have 
\begin{equation}\label{hodge}
h_{i\bar{j}}=\lambda 
g_{i\bar{j}}+g^{\alpha\bar\beta} F_{i\bar j\alpha\bar\beta},
\end{equation}
where $\lambda=\mu-m-1$. Thus  
$\omega_{\mu}>0$ which implies $\omega_{\mu}$ is \ka.

Usually, choosing a normal coordinate system
will simplify the computation greatly. However, in the following computation, 
the use of general coordinates will make the computation easier.

To simplify the computation we first calculate
$\bar\partial_{l}D_{\alpha}D_{i}\Omega$. Using Remark \ref{setconvention} and 
Lemma \ref{h31} we have  
\begin{align}\label{h3rd}
\begin{split}
\bar\partial_{l}D_{\alpha}D_{i}\Omega
=&\bar\partial_{l}(\partial_{\alpha}D_{i}\Omega+K_{\alpha}D_{i}\Omega
-\Gamma^{\gamma}_{i\alpha}D_{\gamma}\Omega)  \\
=&\partial_{\alpha}(\bar\partial_{l}D_{i}\Omega)+(\bar\partial_{
l}K_{\alpha}) D_{i}\Omega+K_{\alpha}\bar\partial_{l}D_{i}\Omega
-(\bar\partial_{l}\Gamma^{\gamma}_{i\alpha})D_{\gamma}\Omega
-\Gamma^{\gamma}_{i\alpha}\bar\partial_{l}D_{\gamma}\Omega 
\\
=& \partial_{\alpha} (g_{i\bar l}\Omega)+g_{\alpha\bar l}D_{i}\Omega
+K_{\alpha}g_{i\bar l}\Omega
-R_{i\bar\tau \alpha\bar l}g^{\gamma\bar\tau}
D_{\gamma}\Omega 
-g^{\gamma\bar\tau}
\frac{\partial g_{i\bar\tau}}{\partial z_{\alpha}}g_{\gamma\bar l}\Omega 
\\
=& g_{i\bar l}D_{\alpha}\Omega+g_{\alpha\bar l}D_{i}\Omega
-R_{i\bar\tau \alpha\bar l}g^{\gamma\bar\tau}D_{\gamma}\Omega\\
=& F_{i\bar\tau \alpha\bar l}g^{\gamma\bar\tau}D_{\gamma}\Omega.
\end{split}
\end{align}

Similarly, we have
\begin{equation}\label{h3rdbar1}
\partial_{k}\bar{D_{\beta}D_{j}\Omega}
=\bar{\bar\partial_{k}D_{\beta}D_{j}\Omega}
=F_{p\bar\beta k\bar j}
g^{p\bar q}\bar{D_{q}\Omega} 
\end{equation}
since $F$ is a curvature like tensor. 
Now because $\Omega$ is holomorphic, we have $\partial_{k}\bar\Omega=
\bar\partial_{l}\Omega=0$. Using Lemma~\ref{h22}, equation ~\eqref{h3rdbar1} 
and the Hodge-Riemann relations,  we know 
\[
(D_\alpha D_i\Omega,\partial_{k}\bar{D_{\beta}D_{j}\Omega})
=F_{p\bar\beta k\bar j}g^{p\bar q}(D_\alpha D_i\Omega,\bar{D_{q}\Omega})=0
\]
which implies 
\begin{align}\label{h1st}
\begin{split}
\frac{\partial h_{i\bar j}}{\partial z_{k}}=&
\lambda\frac{\partial g_{i\bar j}}{\partial z_{k}}+
\frac{(\partial_{k}D_{\alpha}D_{i}\Omega,\bar{D_{\beta}D_{j}\Omega})}
{(\Omega,\bar{\Omega})}g^{\alpha\bar\beta}+
\frac{(D_{\alpha}D_{i}\Omega,\partial_{k}\bar{D_{\beta}D_{j}\Omega})}
{(\Omega,\bar{\Omega})}g^{\alpha\bar\beta} \\
&-\frac{(D_{\alpha}D_{i}\Omega,\bar{D_{\beta}D_{j}\Omega})}
{(\Omega,\bar{\Omega})^{2}}g^{\alpha\bar\beta}(\partial_{k}\Omega,\bar\Omega)
+\frac{(D_{\alpha}D_{i}\Omega,\bar{D_{\beta}D_{j}\Omega})}
{(\Omega,\bar{\Omega})}\frac{\partial g^{\alpha\bar\beta}}{\partial z_{k}}\\
=& \lambda\frac{\partial g_{i\bar j}}{\partial z_{k}}+
\frac{(\partial_{k}D_{\alpha}D_{i}\Omega,\bar{D_{\beta}D_{j}\Omega})}
{(\Omega,\bar{\Omega})}g^{\alpha\bar\beta}+
\frac{(K_{k}D_{\alpha}D_{i}\Omega,\bar{D_{\beta}D_{j}\Omega})}
{(\Omega,\bar{\Omega})}g^{\alpha\bar\beta}\\
&-\frac{(\Gamma_{\alpha k}^{p}D_{p}D_{i}\Omega,\bar{D_{\beta}D_{j}\Omega})}
{(\Omega,\bar{\Omega})}g^{\alpha\bar\beta}\\
=& \lambda\frac{\partial g_{i\bar j}}{\partial z_{k}}+
\frac{(T_{k\alpha i},\bar{D_{\beta}D_{j}\Omega})}
{(\Omega,\bar{\Omega})}g^{\alpha\bar\beta}+\Gamma_{ik}^{p}
\frac{(D_{p}D_{\alpha}\Omega,\bar{D_{\beta}D_{j}\Omega})}
{(\Omega,\bar{\Omega})}g^{\alpha\bar\beta}\\
=& \lambda\frac{\partial g_{i\bar j}}{\partial z_{k}}+
\frac{(T_{k\alpha i},\bar{D_{\beta}D_{j}\Omega})}
{(\Omega,\bar{\Omega})}g^{\alpha\bar\beta}+\Gamma_{ik}^{p}
(h_{p\bar j}-\lambda g_{p\bar j})\\
=& \frac{(T_{k\alpha i},\bar{D_{\beta}D_{j}\Omega})}
{(\Omega,\bar{\Omega})}g^{\alpha\bar\beta}+\Gamma_{ik}^{p}
h_{p\bar j}.
\end{split}
\end{align}
Similarly we have 
\begin{eqnarray}\label{h1stbar}
\frac{\partial h_{i\bar j}}{\partial \bar z_{l}}=
\frac{(D_{\alpha}D_{i}\Omega,\bar{T_{l\beta j}})}
{(\Omega,\bar{\Omega})}g^{\alpha\bar\beta}+\bar{\Gamma_{jl}^{q}}
h_{i\bar q}.
\end{eqnarray}
From \eqref{h1st} and \eqref{h1stbar} we have
\begin{align}\label{h2nd}
\begin{split}
\frac{\partial^{2} h_{i\bar j}}{\partial z_{k}\partial \bar z_{l}}=&
\frac{(\bar\partial_{l}T_{k\alpha i},
\bar{D_{\beta}D_{j}\Omega})}{(\Omega,\bar{\Omega})}g^{\alpha\bar\beta}
+\frac{(T_{k\alpha i},
\bar{\partial_{l}D_{\beta}D_{j}\Omega})}{(\Omega,\bar{\Omega})}
g^{\alpha\bar\beta} 
 -\frac{(T_{k\alpha i}\Omega,\bar{D_{\beta}D_{j}\Omega})}
{(\Omega,\bar{\Omega})^{2}}g^{\alpha\bar\beta}
(\Omega,\bar{\partial_{l}\Omega})\\
&+\frac{(T_{k\alpha i},\bar{D_{\beta}D_{j}\Omega})}
{(\Omega,\bar{\Omega})}\frac{\partial g^{\alpha\bar\beta}}{\partial \bar z_{
l}}
 +(\bar\partial_{l}\Gamma_{ik}^{p})h_{p\bar j}+\Gamma_{ik}^{p}
\frac{\partial h_{p\bar j}}{\partial \bar z_{l}}\\
=& \frac{(\bar\partial_{l}T_{k\alpha i},
\bar{D_{\beta}D_{j}\Omega})}{(\Omega,\bar{\Omega})}g^{\alpha\bar\beta}
+\frac{(T_{k\alpha i},
\bar{\partial_{l}D_{\beta}D_{j}\Omega})}{(\Omega,\bar{\Omega})}
g^{\alpha\bar\beta} 
+\frac{(T_{k\alpha i},\bar{K_{l}D_{\beta}D_{j}\Omega})}
{(\Omega,\bar{\Omega})}g^{\alpha\bar\beta}\\
&-\frac{(T_{k\alpha i},
\bar{\Gamma_{\beta l}^{q}D_{q}D_{j}\Omega})}
{(\Omega,\bar{\Omega})}g^{\alpha\bar\beta}
+R_{i\bar q k\bar l}g^{p\bar q}h_{p\bar j}
+\Gamma_{ik}^{p}\bigg ( \frac{(D_{\alpha}D_{p}\Omega,\bar{T_{l\beta j}})}
{(\Omega,\bar{\Omega})}g^{\alpha\bar\beta}+\bar{\Gamma_{jl}^{q}}h_{p\bar q}
\bigg )\\
=& \frac{(\bar\partial_{l}T_{k\alpha i},
\bar{D_{\beta}D_{j}\Omega})}{(\Omega,\bar{\Omega})}g^{\alpha\bar\beta}
+\frac{(T_{k\alpha i},\bar{T_{l\beta j}})}{(\Omega,\bar{\Omega})}
g^{\alpha\bar\beta} 
+\bar{\Gamma_{jl}^{q}}\frac{(T_{k\alpha i},
\bar{D_{q}D_{\beta}\Omega})}
{(\Omega,\bar{\Omega})}g^{\alpha\bar\beta}\\
&+R_{i\bar q k\bar l}g^{p\bar q}h_{p\bar j}
+\Gamma_{ik}^{p}\bigg ( \frac{(D_{\alpha}D_{p}\Omega,\bar{T_{l\beta j}})}
{(\Omega,\bar{\Omega})}g^{\alpha\bar\beta}+\bar{\Gamma_{jl}^{q}}h_{p\bar q}
\bigg ) .
\end{split}
\end{align}
Since 
\begin{eqnarray*}
\widetilde{R}_{i\bar j k\bar l}=
\frac{\partial^{2} h_{i\bar j}}{\partial z_{k}\partial \bar z_{l}}-
h^{s\bar t}\frac{\partial h_{i\bar t}}{\partial z_{k}}
\frac{\partial h_{s\bar j}}{\partial \bar z_{l}},
\end{eqnarray*}
by \eqref{h2nd}, \eqref{h1st} and \eqref{h1stbar}, we have 
\begin{align}\label{curv10}
\begin{split}
\widetilde{R}_{i\bar j k\bar l}=&
\frac{(\bar\partial_{l}T_{k\alpha i},
\bar{D_{\beta}D_{j}\Omega})}{(\Omega,\bar{\Omega})}g^{\alpha\bar\beta}
+\frac{(T_{k\alpha i},\bar{T_{l\beta j}})}{(\Omega,\bar{\Omega})}
g^{\alpha\bar\beta} 
+\bar{\Gamma_{jl}^{q}}\frac{(T_{k\alpha i},
\bar{D_{q}D_{\beta}\Omega})}
{(\Omega,\bar{\Omega})}g^{\alpha\bar\beta}\\
&+R_{i\bar q k\bar l}g^{p\bar q}h_{p\bar j}
+\Gamma_{ik}^{p}\bigg ( \frac{(D_{\alpha}D_{p}\Omega,\bar{T_{l\beta j}})}
{(\Omega,\bar{\Omega})}g^{\alpha\bar\beta}+\bar{\Gamma_{jl}^{q}}h_{p\bar q}
\bigg )\\
&-h^{s\bar t}\bigg ( \frac{(T_{k\alpha i},\bar{D_{\beta}D_{t}\Omega})}
{(\Omega,\bar{\Omega})}g^{\alpha\bar\beta}+\Gamma_{ik}^{p}h_{p\bar t}
\bigg )
\bigg ( \frac{(D_{\gamma}D_{s}\Omega,\bar{T_{l\tau j}})}
{(\Omega,\bar{\Omega})}g^{\gamma\bar\tau}+\bar{\Gamma_{jl}^{q}}h_{s\bar q}
\bigg )\\
=& \frac{(\bar\partial_{l}T_{k\alpha i},
\bar{D_{\beta}D_{j}\Omega})}{(\Omega,\bar{\Omega})}g^{\alpha\bar\beta}
+\frac{(T_{k\alpha i},\bar{T_{l\beta j}})}{(\Omega,\bar{\Omega})}
g^{\alpha\bar\beta} 
+\bar{\Gamma_{jl}^{q}}\frac{(T_{k\alpha i},
\bar{D_{q}D_{\beta}\Omega})}
{(\Omega,\bar{\Omega})}g^{\alpha\bar\beta}\\
&+R_{i\bar q k\bar l}g^{p\bar q}h_{p\bar j}+\Gamma_{ik}^{p}
\frac{(D_{\alpha}D_{p}\Omega,\bar{T_{l\beta j}})}
{(\Omega,\bar{\Omega})}g^{\alpha\bar\beta}
+\Gamma_{ik}^{p}\bar{\Gamma_{jl}^{q}}h_{p\bar q}\\
&-h^{s\bar t}\frac{(T_{k\alpha i},\bar{D_{\beta}D_{t}\Omega})}
{(\Omega,\bar{\Omega})}g^{\alpha\bar\beta}
\frac{(D_{\gamma}D_{s}\Omega,\bar{T_{l\tau j}})}
{(\Omega,\bar{\Omega})}g^{\gamma\bar\tau}
-h^{s\bar t}\Gamma_{ik}^{p}h_{p\bar t}\bar{\Gamma_{jl}^{q}}h_{s\bar q}\\
&-h^{s\bar t}\frac{(T_{k\alpha i},\bar{D_{\beta}D_{t}\Omega})}
{(\Omega,\bar{\Omega})}g^{\alpha\bar\beta}\bar{\Gamma_{jl}^{q}}h_{s\bar q}
-h^{s\bar t}\frac{(D_{\gamma}D_{s}\Omega,\bar{T_{l\tau j}})}
{(\Omega,\bar{\Omega})}g^{\gamma\bar\tau}\Gamma_{ik}^{p}h_{p\bar t}\\
=& \frac{(T_{k\alpha i},\bar{T_{l\beta j}})}{(\Omega,\bar{\Omega})}
g^{\alpha\bar\beta} 
-h^{s\bar t}\frac{(T_{k\alpha i},\bar{D_{\beta}D_{t}\Omega})}
{(\Omega,\bar{\Omega})}g^{\alpha\bar\beta}
\frac{(D_{\gamma}D_{s}\Omega,\bar{T_{l\tau j}})}
{(\Omega,\bar{\Omega})}g^{\gamma\bar\tau}\\
&+\frac{(\bar\partial_{l}T_{k\alpha i},
\bar{D_{\beta}D_{j}\Omega})}{(\Omega,\bar{\Omega})}g^{\alpha\bar\beta}
+R_{i\bar q k\bar l}g^{p\bar q}h_{p\bar j}.
\end{split}
\end{align}
Using \eqref{h3rd}, Lemma \ref{h31}, Lemma \ref{h22} and the Hodge-Riemann
relations we have 
\begin{align}\label{II}
\begin{split}
&\frac{(\bar\partial_{l}T_{k\alpha i},
\bar{D_{\beta}D_{j}\Omega})}{(\Omega,\bar{\Omega})}g^{\alpha\bar\beta}
\\
=&\frac{(\partial_{k}\bar\partial_{l}D_{\alpha}D_{i}\Omega,
\bar{D_{\beta}D_{j}\Omega})}{(\Omega,\bar{\Omega})}g^{\alpha\bar\beta}
+\frac{((\bar\partial_{l}K_{k})D_{\alpha}D_{i}\Omega,
\bar{D_{\beta}D_{j}\Omega})}{(\Omega,\bar{\Omega})}g^{\alpha\bar\beta}\\
&-\frac{((\bar\partial_{l}\Gamma_{\alpha k}^{p})D_{p}D_{i}\Omega,
\bar{D_{\beta}D_{j}\Omega})}{(\Omega,\bar{\Omega})}g^{\alpha\bar\beta}
-\frac{((\bar\partial_{l}\Gamma_{ik}^{p})D_{p}D_{\alpha}\Omega,
\bar{D_{\beta}D_{j}\Omega})}{(\Omega,\bar{\Omega})}g^{\alpha\bar\beta}
\\
=&F_{\alpha \bar l i\bar q}F_{p\bar\beta k\bar j}
g^{\alpha\bar\beta}g^{p\bar q}
+F_{\alpha\bar\beta i\bar j}g^{\alpha\bar\beta}g_{k\bar l}\\
&-R_{\alpha\bar q k\bar l}F_{p\bar\beta i\bar j}
g^{\alpha\bar\beta}g^{p\bar q}
-R_{i\bar q k\bar l}F_{p\bar\beta\alpha\bar j}
g^{\alpha\bar\beta}g^{p\bar q}.
\end{split}
\end{align}
By \eqref{hodge} and the Strominger formula~\eqref{stro1}, the above expression is
\begin{equation}\label{II1}
F_{i\bar q \alpha\bar l}F_{p\bar j k\bar\beta}
g^{\alpha\bar\beta}g^{p\bar q}
+F_{\alpha\bar q k\bar l}F_{i\bar j p\bar\beta}
g^{\alpha\bar\beta}g^{p\bar q}
+\lambda (g_{i\bar j}g_{k\bar l}+g_{i\bar l}g_{k\bar j})
-(\lambda+1)F_{i\bar j k\bar l}.
\end{equation}
Using the Hodge-Riemann relations we have 
\begin{align}\label{I1}
\begin{split}
\frac{(T_{k\alpha i},\bar{T_{l\beta j}})}{(\Omega,\bar{\Omega})}
g^{\alpha\bar\beta} =
\frac{(D_{k}D_{\alpha}D_{i}\Omega,
\bar{D_{l}D_{\beta}D_{j}\Omega})}{(\Omega,\bar\Omega)}g^{\alpha\bar\beta}
+\frac{(E_{k\alpha i},
\bar{E_{l\beta j}})}{(\Omega,\bar\Omega)}g^{\alpha\bar\beta}
\end{split}
\end{align}
and 
\begin{align}\label{I2}
\begin{split}
 h^{s\bar t}\frac{(T_{k\alpha i},\bar{D_{\beta}D_{t}\Omega})}
{(\Omega,\bar\Omega)}\frac{(D_{\gamma}D_{s}\Omega,
\bar{T_{l\tau j}})}{(\Omega,\bar\Omega)}g^{\alpha\bar\beta}g^{\gamma\bar\tau}
= h^{s\bar t}\frac{(E_{k\alpha i},\bar{D_{\beta}D_{t}\Omega})}
{(\Omega,\bar\Omega)}\frac{(D_{\gamma}D_{s}\Omega,
\bar{E_{l\tau j}})}{(\Omega,\bar\Omega)}g^{\alpha\bar\beta}g^{\gamma\bar\tau}.
\end{split}
\end{align}
Theorem~\ref{curvmu} follows from
\eqref{II1} \eqref{I1}, and \eqref{I2}.

\qed


\bibliographystyle{abbrv} 
\bibliography{bib}

\begin{thebibliography}{10}

\bibitem{BG}
R.~Bryant and P.~Griffiths.
\newblock Some observations on the infinitesimal period relations for regular
  threefolds with trivial canonical bundle.
\newblock In M.~Artin and J.~Tate, editors, {\em Arithmetic and Geometry},
  pages 77--85. Boston, Birkha\"user, 1983.

\bibitem{COGP}
P.~Candelas, X.~C. de~la Ossa, P.~S. Green, and L.~Parkes.
\newblock A pair of {C}alabi-{Y}au manifolds as an exactly soluble
  superconformal theory.
\newblock In S.Y.Yau, editor, {\em Essays in mirror manifolds}, pages 31--95.
  International Press, 1992.

\bibitem{Fd}
R.~Friedman.
\newblock {On threefolds with trivial canonical bundle}.
\newblock In S.-T. Yau, editor, {\em Complex geometry and Lie theory (Sundance,
  UT, 1989)}, pages 103--134. International Press, 1992.

\bibitem{Gr}
P.~Griffiths, editor.
\newblock {\em Topics in Transcendental Algebraic Geometry}, volume 106 of {\em
  Ann. Math Studies}.
\newblock Princeton University Press, 1984.

\bibitem{GS}
P.~Griffiths and W.~Schmid.
\newblock {Locally homogeneous complex manifolds}.
\newblock {\em Acta Math}, 123:253--302, 1969.

\bibitem{Gr1}
P.~A. Griffiths.
\newblock Periods of integrals on algebraic manifolds. {I}. {C}onstruction and
  properties of the modular varieties.
\newblock {\em Amer. J. Math.}, 90:568--626, 1968.

\bibitem{Gr2}
P.~A. Griffiths.
\newblock Periods of integrals on algebraic manifolds. {I}{I}. {L}ocal study of
  the period mapping.
\newblock {\em Amer. J. Math.}, 90:805--865, 1968.

\bibitem{Hi}
H.~Hironaka.
\newblock Resolution of singularities of an algebraic variety over a field of
  characteristic 0.
\newblock {\em Ann. of Math.}, 79:109--326, 1964.

\bibitem{JY1}
J.~Jost and S.-T. Yau.
\newblock {Harmonic mappings and algebraic varieties over function fields}.
\newblock {\em Amer. J. Math}, 115(6):1197--1227, 1993.

\bibitem{LTYZ}
K.~Liu, A.~Todorov, S.~Yau, and K.~Zuo.
\newblock Shafarevich type theorem of moduli space of {C}alabi-{Y}au manifolds.
\newblock preprint.

\bibitem{Lu3}
Z.~Lu.
\newblock On the geometry of classifying spaces and horizontal slices.
\newblock {\em Amer. J. Math.}, 121:177--198, 1999.

\bibitem{Lu5}
Z.~Lu.
\newblock On the {H}odge metric of the universal deformation space of
  {Calabi-Yau} threefolds.
\newblock {\em J. Geom. Anal.}, 11(1):103--118, 2001.

\bibitem{LS-2}
Z.~Lu and X.~Sun.
\newblock The {W}eil-{P}etersson volume of the moduli space of {C}alabi-{Y}au
  manifolds.
\newblock preprint, 2002.

\bibitem{MZograf}
Y.~Manin and P.~Zograf.
\newblock Invertible cohomological field theories and {W}eil-{P}etersson
  volumes.
\newblock AG/9902051.

\bibitem{McM}
C.~T. McMullen.
\newblock The moduli space of {R}iemann surfaces is {K}\"ahler hyperbolic.
\newblock {\em Ann. of Math.}, 151:327--357, 2000.

\bibitem{Mok1}
N.~Mok.
\newblock Compactification of complete {K}\"ahler surface of finite volume
  satisfying certain conditions.
\newblock {\em Ann. Math}, 128:383--425, 1989.

\bibitem{Mumford1}
D.~Mumford.
\newblock Hirzebruch's proportionality theorem in the non-compact case.
\newblock {\em Invent. Math}, 42:239--272, 1977.

\bibitem{KMZ}
K.~R., Y.~Manin, and D.~Zagier.
\newblock Higher {W}eil-{P}etersson volumes of moduli spaces of stable
  $n$-pinted curvatures.
\newblock {\em Comm. Math. Phys}, 181:763--787, 1996.

\bibitem{Royden}
H.~L. Royden.
\newblock The {A}hlfors-{S}chwarz lemma in several complex variables.
\newblock {\em Comment. Math. Helv.}, 55(4):547--558, 1980.

\bibitem{Schmid}
W.~Schmid.
\newblock Variation of {H}odge structure: the singularities of the period
  mapping.
\newblock {\em Inventiones Math.}, 22:211--319, 1973.

\bibitem{Sc1}
G.~Schumacher.
\newblock The curvature of the {P}etersson-{W}eil metric on the moduli space of
  {K\"ahler-Einstein} manifolds.
\newblock In V.~Ancona and A.~Silva, editors, {\em Complex Analysis and
  Geometry}, pages 339--345, 1993.

\bibitem{Sc2}
G.~Schumacher and S.~Trapani.
\newblock Estimates of {W}eil-{P}etersson volumes via effective divisors.
\newblock {\em Comm. Math. Phys.}, 222:1--7, 2001.

\bibitem{Siu2}
Y.~T. Siu.
\newblock Curvature of the {Weil-Petersson} metric in the moduli spaces of
  compact {K\"ahler-Einstein} manifolds of negative first chern class.
\newblock In P.-M. Wong and A.~Howard, editors, {\em Complex Analysis, Papers
  in Honour of Wilhelm Stoll}. Vieweg, Braunschweig, 1986.

\bibitem{S}
A.~Strominger.
\newblock {Special Geometry}.
\newblock {\em Comm. Math. Phy.}, 133:163--180, 1990.

\bibitem{T1}
G.~Tian.
\newblock Smoothness of the universal deformation space of compact {Calabi-Yau}
  manifolds and its {Peterson-Weil} metric.
\newblock In S.-T. Yau, editor, {\em Mathematical aspects of string theory},
  volume~1, pages 629--646. World Scientific, 1987.

\bibitem{To}
A.~N. Todorov.
\newblock {The Weil-Petersson geometry of the moduli space of ${\rm SU}(n\geq
  3)$ (Calabi-Yau) manifolds. I}.
\newblock {\em Comm. Math. Phys}, 126(2):325--346, 1989.

\bibitem{V2}
E.~Viehweg.
\newblock {\em Quasi-projective moduli for polarized manifolds}.
\newblock Ergebnisse der Mathematik und ihrer Grenzgebiete. Springer-Verlag,
  1991.

\bibitem{Wang1}
C.-L. Wang.
\newblock Curvature properties of the {C}alabi-{Y}au moduli.
\newblock to appear in Documenta Mathematica.

\bibitem{Wong}
C.-L. Wang.
\newblock On the incompleteness of the {Weil-Petersson} metrics along
  degenerations of {Calabi-Yau} manifolds.
\newblock {\em Math. Res. Lett.}, 1:157--171, 1997.

\bibitem{Wolpert}
S.~Wolpert.
\newblock Chern forms and the {R}iemann tensor for the moduli space of curves.
\newblock {\em Invent. Math}, 85:119--145, 1986.

\bibitem{Y3}
S.-T. Yau.
\newblock A general {S}chwarz lemma for {K}\"ahler manifolds.
\newblock {\em Amer. J. Math.}, 100:197--203, 1978.

\bibitem{Zograf3}
P.~Zograf.
\newblock Weil-{P}etersson volumes of low genus moduli spaces of curves and
  genus expansion in two dimensional gravity.
\newblock AG/9811026.

\bibitem{Zograf1}
P.~Zograf.
\newblock The {W}eil-{P}etersson volume of the moduli space of punctured
  spheres.
\newblock In B\"odigheimer, editor, {\em Mapping class groups and moduli spaces
  of {R}iemann surfaces}, volume 150 of {\em Contemp. Math}, pages 367--372.
  American Math Soc., 1993.

\bibitem{Zograf2}
P.~Zograf.
\newblock Weil-{P}etersson volumes of low genus moduli spaces.
\newblock {\em Funct. Anal. Appl.}, 32:78--81, 1998.

\end{thebibliography}

\end{document}